# On explicit lifts of cusp forms from $\mathrm{GL}_m$ to classical groups

By David Ginzburg, Stephen Rallis, and David Soudry*

## Introduction

In this paper, we begin the study of poles of partial $L$-functions $L^S(\sigma \otimes \tau, s)$, where $\sigma \otimes \tau$ is an irreducible, automorphic, cuspidal, generic (i.e. with nontrivial Whittaker coefficient) representation of $G_{\mathbb{A}} \times \mathrm{GL}_m(\mathbb{A})$. $G$ is a split classical group and $\mathbb{A}$ is the adele ring of a number field $F$. We also consider $\widetilde{\mathrm{Sp}}_{2n}(\mathbb{A}) \times \mathrm{GL}_m(\mathbb{A})$, where $\sim$ denotes the metaplectic cover.

Examining $L^S(\sigma \otimes \tau, s)$ through the corresponding Rankin-Selberg, or Shimura-type integrals ([G-PS-R], [G-R-S3], [G], [So]), we find that the global integral contains, in its integrand, a certain normalized Eisenstein series which is responsible for the poles. For example, if $G = \mathrm{SO}_{2k+1}$, then the Eisenstein series is on the adele points of split $\mathrm{SO}_{2m}$, induced from the Siegel parabolic subgroup and $\tau \otimes |\det \cdot|^{s-1/2}$. If $G = \widetilde{\mathrm{Sp}}_{2k}$ (this is a convenient abuse of notation), then the Eisenstein series is on $\mathrm{Sp}_{2m}(\mathbb{A})$, induced from the Siegel parabolic subgroup and $\tau \otimes |\det \cdot|^{s-1/2}$. The constant term (along the "Siegel radical") of such a normalized Eisenstein series involves one of the $L$-functions $L^S(\tau, \Lambda^2, 2s-1)$ or $L^S(\tau, \mathrm{Sym}^2, 2s-1)$. So, up to problems of normalization of intertwining operators, the only pole we expect, for $\mathrm{Re}(s) > 1/2$, is at $s = 1$ and then $\tau$ should be self-dual. (See [J-S1], [B-G].) Thus, let us assume that $\tau$ is self-dual. By [J-S1], we know that $L^S(\tau \otimes \tau, s)$ has a simple pole at $s = 1$. Since

$$L^S(\tau \otimes \tau, s) = L^S(\tau, \Lambda^2, s) L^S(\tau, \mathrm{Sym}^2, s) \tag{0.1}$$

and since each factor is nonzero at $s = 1$, it follows that exactly one of the $L$-functions $L^S(\tau, \Lambda^2, s)$ and $L^S(\tau, \mathrm{Sym}^2, s)$ has a simple pole at $s = 1$. We recall that if $m$ is odd, then $L^S(\tau, \Lambda^2, s)$ is entire ([J-S1]).

Assume, for example, that $m = 2n$ and $L^S(\tau, \Lambda^2, s)$ has a pole at $s = 1$. Then the above-mentioned Eisenstein series (denote it, for this introduction, by $E_{\tau,s}(h)$) on $\mathrm{SO}_{4n}(\mathbb{A})$ has a simple pole at $s = 1$. (See Proposition 1 and Remark 2.) Let $G = \mathrm{SO}_{2k+1}$ and $\sigma$ be as above. Langlands conjectures, predict

*This research was supported by the Basic Research Foundation administered by the Israel Academy of Sciences and Humanities.

the existence of a functorial lift of $\sigma$, denote it $\pi$, on $\mathrm{GL}_{2k}(\mathbb{A})$. $\pi$ is an irreducible automorphic representation of $\mathrm{GL}_{2k}(\mathbb{A})$. In particular, $L^S(\sigma\otimes\tau,s) = L^S(\pi\otimes\tau,s)$. Thus, if $k<n$, $L^S(\sigma\otimes\tau,s)$ should be holomorphic at $s=1$. Similarly, if $L^S(\tau,\Lambda^2,s)$ has a pole at $s=1$, $(m=2n)$, and $L^S(\tau,st,1/2) \neq 0$, (partial standard $L$-function) then the above-mentioned Eisenstein series (denote it by $E_{\tau,s}(h)$) on $\widetilde{\mathrm{Sp}}_{4n}(\mathbb{A})$ has a simple pole at $s=1$ (see Proposition 1). Let $G=\widetilde{\mathrm{Sp}}_{2k}$ and $\sigma$ as above. Again, one expects the existence of a "functorial lift" of $\sigma$, denote it $\pi$, on $\mathrm{GL}_{2k}(\mathbb{A})$ (such a lift is not canonic, it depends on a choice of a nontrivial character $\psi$ of $F\backslash\mathbb{A}$). One of our main results (Theorem 14) states

THEOREM A. *Let $\sigma\otimes\tau$ be an irreducible, automorphic, cuspidal, generic representation of $\widetilde{Sp}_{2k}(\mathbb{A})\times\mathrm{GL}_{2n}(\mathbb{A})$, such that $1\le k<n$. Assume that*

$$L^S(\tau,st,1/2)\neq 0 \tag{0.2}$$

*and*

$$L^S(\tau,\Lambda^2,s) \quad \text{has a pole at} \quad s=1\ .$$

*Then $L^S_\psi(\sigma\otimes\tau,s)$ is holomorphic at $s=1$.*

The definition of the partial standard $L$-function of $\sigma\otimes\tau$ depends on a choice of $\psi$; see [G-R-S3]. The proof of the theorem involves new ideas and results which we regard as the main contribution of this work.

In general terms, we actually construct a candidate $\sigma=\sigma(\tau)=\sigma_\psi(\tau)$, on $\widetilde{\mathrm{Sp}}_{2n}(\mathbb{A})$, which should lift functorially to $\tau$. Moreover, $\sigma(\tau)$ fits into a fascinating "tower" of automorphic cuspidal modules $\sigma_k(\tau)$ of $\widetilde{\mathrm{Sp}}_{2k}(\mathbb{A})$, $0\le k\le 2n-1$. (Define $\widetilde{\mathrm{Sp}}_0(\mathbb{A})=\{1\}$.) The "tower property" is the same as in the theory of theta series liftings of dual pairs [R]. Thus, for the first index $\ell_n$, such that $\sigma_{\ell_n}(\tau)$ is nontrivial, $\sigma_{\ell_n}(\tau)$ is cuspidal, and for higher indices, $\sigma_k(\tau)$ is noncuspidal.

Let us give some more details. We first consider the example where $\sigma$ is a generic representation of $\mathrm{SO}_{2k+1}(\mathbb{A})$ and $\tau$ is on $\mathrm{GL}_{2n}(\mathbb{A})$, such that $L^S(\tau,\Lambda^2,s)$ has a pole at $s=1$. The Rankin-Selberg integral for $L^S(\sigma\otimes\tau,s)$, when $k<2n$, has the form

$$\int\limits_{\mathrm{SO}_{2k+1}(F)\backslash\mathrm{SO}_{2k+1}(\mathbb{A})} \varphi(g)E^{\psi_k}_{\tau,s}(g)dg \tag{0.3}$$

where $E^{\psi_k}_{\tau,s}$ is a certain Fourier coefficient of the Eisenstein series (on $\mathrm{SO}_{4n}(\mathbb{A})$) along some unipotent subgroup $N_k$ of $\mathrm{SO}_{4n}$, and with respect to a character $\psi_k$ of $N_k(F)\backslash N_k(\mathbb{A})$, stabilized by $\mathrm{SO}_{2k+1}(\mathbb{A})$, for some corresponding embedding of $\mathrm{SO}_{2k+1}$ inside $\mathrm{SO}_{4n}$. $\varphi$ is a cusp form in the space of $\sigma$. The integral (0.3), in case $k = n$, suggests that if $L^S(\sigma \otimes \tau, s)$ has a pole at $s = 1$, then $\sigma$ pairs with the representation $\sigma(\tau)$, of $\mathrm{SO}_{2n+1}(\mathbb{A})$, on the space spanned by

$$\left[ \mathrm{Res}_{s=1} \overline{E_{\tau,s}(h)} \right]^{\psi_n} \Bigg|_{\mathrm{SO}_{2n+1}(\mathbb{A})} .$$

Now, when we try to compute the various constant terms of $\sigma(\tau)$ along the unipotent radicals of parabolic subgroups of $\mathrm{SO}_{2n+1}$, we find out that these are expressed in terms of

$$\left[ \mathrm{Res}_{s=1} \overline{E_{\tau,s}(h)} \right]^{\psi_\ell} \Bigg|_{\mathrm{SO}_{2\ell+1}(\mathbb{A})}$$

for $\ell < n$. This immediately reveals the tower property. Indeed, we define for any $0 \le k < 2n$, the representation $\sigma_k(\tau)$ of $\mathrm{SO}_{2k+1}(\mathbb{A})$, acting by right translations on the space of

$$\left[ \mathrm{Res}_{s=1} \overline{E_{\tau,s}(h)} \right]^{\psi_k} \Bigg|_{\mathrm{SO}_{2k+1}(\mathbb{A})} .$$

This leads us to another main result of this paper.

THEOREM B. *Let $\tau$ be an irreducible, automorphic, cuspidal representation of $\mathrm{GL}_{2n}(\mathbb{A})$. Assume that $L^S(\tau, \Lambda^2, s)$ has a pole at $s = 1$. Then the representations $\{\sigma_k(\tau)\}_{k=0}^{2n-1}$ have the tower property, i.e. for the first index $\ell_n$, such that $\sigma_{\ell_n}(\tau) \neq 0$, $\sigma_{\ell_n}(\tau)$ is cuspidal and for $k > \ell_n$, $\sigma_k(\tau)$ is noncuspidal.*

It is easy to see that $1 \le \ell_n \le 2n-1$ and that $\sigma_{\ell_n}(\tau)$ is generic. Our main conjecture for this case is

CONJECTURE.
1) $\ell_n = n$.
2) *$\sigma_n(\tau)$ is irreducible and lifts functorially from $\mathrm{SO}_{2n+1}(\mathbb{A})$ to $\tau$.*

*Remark.* The conjecture implies that $\sigma_n(\tau)$ is a generic member of the $L$-packet which lifts to $\tau$.

We define similar towers, prove Theorem B and state the above conjecture in case conditions (0.2) hold, and also in case $L^S(\tau, \mathrm{Sym}^2, s)$ has a pole at $s = 1$ (so here $\tau$ is on $\mathrm{GL}_{2n}(\mathbb{A})$ or on $\mathrm{GL}_{2n+1}(\mathbb{A})$). In each case, we use a corresponding global integral. For example, if conditions (0.2) hold, then the global integral is of Shimura type, and we construct the representations $\sigma_k(\tau)$ of $\widetilde{\mathrm{Sp}}_{2k}(\mathbb{A})$ using a sequence of Fourier-Jacobi coefficients of $\mathrm{Res}_{s=1} E_{\tau,s}(h)$ (the Eisenstein series is on $\mathrm{Sp}_{4n}(\mathbb{A})$). In this case, we make one step towards the conjecture and prove

Theorem C. *Let $\tau$ be an irreducible, automorphic, cuspidal representation of* $\mathrm{GL}_{2n}(\mathbb{A})$, *such that* $L^S(\tau, \Lambda^2, s)$ *has a pole at* $s = 1$ *and* $L^S(\tau, st, 1/2) \neq 0$. *Then*

$$\sigma_k(\tau) = 0 \ , \quad \text{for} \quad k < n \ ;$$

*that is* $\ell_n \geq n$.

The proof of Theorem C is based on the following two key observations. The first is

Theorem D. *Under the assumptions of Theorem* C, $\mathrm{Res}_{s=1} E_{\tau,s}(h)$ *has a nontrivial period along the subgroup* $\mathrm{Sp}_{2n}(\mathbb{A}) \times \mathrm{Sp}_{2n}(\mathbb{A})$ (*direct sum embedding in* $\mathrm{Sp}_{4n}(\mathbb{A})$).

The second observation is that for an irreducible, admissible representation of $\mathrm{Sp}_{4n}(F)$, where $F$ is $p$-adic, the existence of (nontrivial) $\mathrm{Sp}_{2n}(F) \times \mathrm{Sp}_{2n}(F)$ invariant functionals negates the existence of the Fourier-Jacobi models which enter the definition of $\sigma_k(\tau)$, provided $k < n$.

*Note added in proof.* Since the time of writing this paper, we have proved large parts of the above conjecture. We can now prove that $\ell_n = n$ in all four cases dealt with here. Thus, for $\tau$ as above, $\sigma_n(\tau)$ is a nonzero cuspidal module on $G_{\mathbb{A}}$ ($G = \widetilde{\mathrm{Sp}}_{2n}$, $\mathrm{SO}_{2n+1}, \mathrm{SO}_{2n}, \mathrm{Sp}_{2n}$ respectively). The proof for $G = \widetilde{\mathrm{Sp}}_{2n}$ appears in [G-R-S5], where we also give the analogous local theory, and construct for an irreducible supercuspidal, self-dual representation $\tau$ of $\mathrm{GL}_{2n}(k)$, $k$ – a $p$-adic field, such that its local exterior square $L$-function has a pole at $s = 0$, an irreducible, supercuspidal, generic (with respect to a prechosen character $\psi$) representation $\sigma = \sigma_\psi(\tau)$, such that $\gamma(\sigma \otimes \tau, s, \psi)$ has a pole at $s = 1$. $\gamma(\sigma \otimes \tau, s, \psi)$ is the local gamma factor attached to $\sigma \otimes \tau$, by the theory of Shimura type integrals considered here. In case $G = \widetilde{\mathrm{Sp}}_{2n}$, we prove in [G-R-S6] that each irreducible summand of $\sigma_n(\tau)$ has indeed the unramified parameters determined by $\tau$, at almost all places (once we fix $\psi$). Moreover, $\sigma_n(\tau)$ is the direct sum of all irreducible, cuspidal, $\psi$-generic representations of $\widetilde{\mathrm{Sp}}_{2n}(\mathbb{A})$, which $\psi$-lift to $\tau$. Thus, if we have the irreducibility of $\sigma_n(\tau)$, which is part of the above conjecture, then $\sigma_n(\tau)$ will be the unique, $\psi$-generic member of the "$\psi$-$L$-packet determined by $\tau$." The proof of the unramified correspondence and of the fact that $\ell_n = n$ in the remaining cases will appear in a paper which is now in preparation.

The results of this paper and those just mentioned were announced in [G-R-S2]; see also [G-R-S1].

Finally, this paper is organized as follows. In Section 1, we prove Theorem D. In Section 2 we define the towers of representations $\sigma_k(\tau)$ and prove Theorem B for all cases $(\mathrm{SO}_{2k+1}, \widetilde{\mathrm{Sp}}_{2k}, \mathrm{SO}_{2k}, \mathrm{Sp}_{2k})$. In Section 3 we prove Theorem C. Theorem A then follows as a corollary.

*Acknowledgement.* We thank Jim Cogdell for helpful discussions and for his interest.

## 1. $\mathrm{Sp}_{2n} \times \mathrm{Sp}_{2n}$ — Periods of residues of Eisenstein series

1. *The Eisenstein series of study and its pole at $s = 1$.* Let $G = \mathrm{Sp}_{4n}$, considered as an algebraic group defined over a number field $F$. Let $P = MU$ be the Siegel parabolic subgroup of $G$. Consider an irreducible, automorphic, cuspidal representation $\tau$ of $\mathrm{GL}_{2n}(\mathbb{A})$. Assume it is self-dual. Regard $\tau$ as a representation of $M_{\mathbb{A}}$ as well. Let $\phi$ be an element of $\mathrm{Ind}_{P_{\mathbb{A}}}^{G_{\mathbb{A}}} \tau$; i.e. $\phi$ is a smooth function on $G_{\mathbb{A}}$, taking values in the space of $\tau$, such that

$$\phi(mug; r) = \delta_P^{1/2}(m)\phi(g; rm) ,$$

for $m \in M_{\mathbb{A}}$, $u \in U_{\mathbb{A}}$, $g \in G_{\mathbb{A}}$, $r \in \mathrm{GL}_{2n}(\mathbb{A})$. We realize $\phi$ as a complex function on $G_{\mathbb{A}} \times \mathrm{GL}_{2n}(\mathbb{A})$, such that $r \mapsto \phi(g; r)$ is a cusp form in the space of $\tau$. Assume that $\phi$ is right $K$-finite, where $K$ is the standard maximal compact subgroup of $G_{\mathbb{A}}$. Let, for $s \in \mathbb{C}$,

$$\begin{aligned} \varphi^{\phi}_{\tau,s}(g; m) &= H(g)^{s-1/2}\phi(g; m) , \quad g \in G_{\mathbb{A}}, m \in \mathrm{GL}_{2n}(\mathbb{A}) \\ f^{\phi}_{\tau,s(g)} &= \varphi^{\phi}_{\tau,s}(g; 1), \end{aligned}$$

where if the Iwasawa decomposition of $g$ is $\widehat{a}uk$, $a \in \mathrm{GL}_{2n}(\mathbb{A})$, $u \in U_{\mathbb{A}}, k \in K$, then $H(g) = |\det a|$. Now consider the corresponding Eisenstein series

$$E\big(g, f^{\phi}_{\tau,s}\big) = \sum_{\gamma \in P_F \backslash G_F} f^{\phi}_{\tau,s}(\gamma g) .$$

The series converges absolutely for $\mathrm{Re}(s) > n + 1$, and admits a meromorphic continuation to the whole plane. It has a finite number of poles in the half-plane $\mathrm{Re}(s) \geq \frac{1}{2}$ and they are all simple [M-W, IV.1.11]. Consider the constant term along $U$,

$$\begin{aligned} E^U(g, f^{\phi}_{\tau,s}) &= \int_{U_F \backslash U_{\mathbb{A}}} E(ug, f^{\phi}_{\tau,s}) du \\ &= f^{\phi}_{\tau,s}(g) + M(s) f^{\phi}_{\tau,s}(g) , \end{aligned} \tag{1.1}$$

where $M(s)$ is the intertwining operator, given, for $\mathrm{Re}(s) > n + 1$, by the

convergent integral

$$M(s)f^{\phi}_{\tau,s}(g) = \int\limits_{U_{\mathbb{A}}} f^{\phi}_{\tau,s}(w^{-1}ug)du \ ,$$

for $w = \begin{pmatrix} & I_{2n} \\ -I_{2n} & \end{pmatrix}$. Later on, it will be convenient to consider the intertwining operator as evaluation at $m = 1$ of

$$M(s)\big(\varphi^{\phi}_{\tau,s}\big)(g;m) = \int\limits_{U_F\backslash U_{\mathbb{A}}} \varphi^{\phi}_{\tau,s}(w^{-1}ug;m)du$$

which converges for $\mathrm{Re}(s) > n+1$. These operators are decomposable in the following sense. Fix realizations $V_{\tau_\nu}$ of the local representations $\tau_\nu$, and fix an isomorphism $\ell : \otimes V_{\tau_\nu} \longrightarrow V_\tau$, where $V_\tau$ is the space of $\tau$. Now write $\phi(g)(x)$ instead of $\phi(g;x)$, for $g \in G_{\mathbb{A}}$, $x \in \mathrm{GL}_{2n}(\mathbb{A})$, so that $\phi(g) \in V_\tau$. For each place $\nu$, let $\phi_\nu$ be an element of $\mathrm{Ind}^{G_\nu}_{P_\nu} \tau_\nu$, such that for almost all $\nu$, $\phi_\nu = \phi^0_\nu$ unramified and $\phi^0_\nu(1) = \xi^0_\nu$ – a prechosen unramified vector in $V_{\tau_\nu}$. Assume that $\phi$ corresponds to $\otimes\phi_\nu$, that is $\ell(\otimes\phi_\nu(g_\nu)) = \phi(g)$. For such $\phi$, we have,

$$M(s)f^{\phi}_{\tau,s}(g) = \ell\Big(\otimes M_\nu(s)\varphi^{\phi_\nu}_{\tau_\nu,s}(g_\nu)\Big)(1) \ , \tag{1.2}$$

where

$$\varphi^{\phi}_{\tau_\nu,s}(g_\nu) = H(g_\nu)^{s-1/2}\phi_\nu(g_\nu)$$

and $M_\nu(s)$ is the vector-valued intertwining operator that is the meromorphic continuation of the absolutely convergent integral given, for $\mathrm{Re}(s) > n+1$, by

$$M_\nu(s)\varphi^{\phi_\nu}_{\tau_\nu,s}(g_\nu) = \int\limits_{U_\nu} \varphi^{\phi_\nu}_{\tau_\nu,s}(w^{-1}ug_\nu)du \ .$$

It is well known that for $\phi_\nu = \phi^0_\nu$

$$M_\nu(s)\varphi^{\phi^0_\nu}_{\tau_\nu,s} = \frac{L(\tau_\nu, st, s-\frac{1}{2})L(\tau_\nu, \Lambda^2, 2s-1)}{L(\tau_\nu, st, s+\frac{1}{2})L(\tau_\nu\Lambda^2, 2s)}\varphi^{\phi^0_\nu}_{\tau_\nu,1-s} \ . \tag{1.3}$$

(Recall that $\widehat{\tau} \cong \tau$.) The factors in (1.3) are respectively the standard and the exterior square local $L$-functions attached to $\tau_\nu$. Let $S$ be a finite set of places outside which $\phi_\nu = \phi^0_\nu$ and $g_\nu \in K_\nu$-the standard maximal compact subgroup of $G_\nu$. Then in (1.2)

(1.4)

$$M(s)f^{\phi}_{\tau,s}(g) = \frac{L^S(\tau, st, s-\frac{1}{2})L^S(\tau, \Lambda^2, 2s-1)}{L^S(\tau, st, s+\frac{1}{2})L^S(\tau, \Lambda^2, 2s)}\ell\Big(\otimes_{\nu\in S} M_\nu(s)\varphi^{\phi_\nu}_{\tau_\nu,s}(g_\nu)\otimes(\xi^0)^S\Big)(1) \ ,$$

where $(\xi^0)^S = \otimes_{\nu\notin S}\xi^0_\nu$.

Proposition 1. *Let $S$ be a finite set of places, including those at infinity, outside which $\tau_\nu$ is unramified. Assume that $L^S(\tau, st, \frac{1}{2}) \neq 0$ and $L^S(\tau, \Lambda^2, s)$ has a pole (simple) at $s = 1$. Then the Eisenstein series $E(g, f^\phi_{\tau,s})$ has a pole (simple) at $s = 1$.*

*Proof.* It is enough to show that the constant term $E^U(g, f^\phi_{\tau,s})$ has a pole at $s = 1$, and, by (1.1), it is enough to show this for $M(s)f^\phi_{\tau,s}(g)$. The $L$-functions $L^S(\tau, st, s+\frac{1}{2})$ and $L(\tau, \Lambda^2, 2s)$, are holomorphic at $s = 1$. (Actually $L^S(\tau, st, s+\frac{1}{2})$ and $L^S(\tau, \Lambda^2, 2s)$ are nonzero at $s = 1$. Indeed, form [J-S2, Th. 5.3], $L^S(\tau, st, z)$ is nonzero for $\mathrm{Re}(z) > 1$. This is also true, by the same theorem, for $L^S(\tau \times \tau, z) = L^S(\tau, \Lambda^2, z)L^S(\tau, \mathrm{sym}^2, z)$. From [J-S1] and [B-G], $L^S(\tau, \Lambda^2, z)$ and $L^S(\tau, \mathrm{sym}^2, z)$ are holomorphic at $z = 2$, and since their product is nonzero at $z = 2$, it follows that each of these $L$-function is holomorphic and nonzero at $z = 2$.) Thus, (1.4) will produce a pole of $M(s)f^\phi_{\tau,s}(g)$ at $s = 1$, provided we can choose data such that $\ell\Big(\otimes_{\nu\in S} M_\nu(s)\varphi^{\phi_\nu}_{\tau_\nu,s}(g_\nu)\otimes(\xi^0)^S\Big)(1)$ is nonzero at $s = 1$. The last expression is holomorphic at a neighbourhood of $s = 1$, since otherwise (1.4) will produce a high-order pole for $M(s)f^\phi_{\tau,s}(g)$, and this implies that $E(g, f^\phi_{\tau,s})$ has a high-order pole at $s = 1$, contradicting the simplicity of the poles of $E(g, f^\phi_{\tau,s})$, for $\mathrm{Re}(s) > 1/2$. It is enough to consider $g = 1$. Let $\nu \in S$ be finite. Choose $\phi_\nu$ to have compact support, modulo $P_\nu$, in the open cell $P_\nu w^{-1}U_\nu$, and such that the function on $U_\nu$, $u \mapsto \phi_\nu(w^{-1}u)$ is the characteristic function of a small neighbourhood in $U_\nu$ of 1. Choose $\xi_\nu \in V_{\tau_\nu}$, so that $\phi_\nu(w^{-1}) = \xi_\nu$. Then $M_\nu(s)\varphi^{\phi_\nu}_{\tau_\nu,s}(1) = c_\nu\xi_\nu$, where $c_\nu$ is the measure (in $U_\nu$) of the neighbourhood above. Let $\nu$ be archimedean. Consider, as above, $\phi_\nu$ which is compactly supported, modulo $P_\nu$, in the open cell $P_\nu w^{-1}U_\nu$, and such that the function $u \mapsto \phi_\nu(w^{-1}u)$ has the form $b_\nu(u)\xi_\nu$, where $b_\nu$ is a Schwarz-Bruhat function on $U_\nu$ and $\xi_\nu \in V_{\tau_\nu}$. $\varphi^{\phi_\nu}_{\tau_\nu,s}$ is a smooth section. Let

$$a(s) = \prod_{\nu \text{ archimedean}} \int_{U_\nu} H_\nu^{s-1/2}(w^{-1}u)b_\nu(u)du\ .$$

Clearly $a(s)$ is holomorphic, and functions $b_\nu$ can be chosen so that $a(1) \neq 0$. For this data,

$$M(s)f^\phi_{\tau,s}(1) = a(s)\frac{L^S(\tau, st, s-\frac{1}{2})L^S(\tau, \Lambda^2, 2s-1)}{L^S(\tau, st, s+\frac{1}{2})L^S(\tau, \Lambda^2, 2s)}\ell(\otimes\xi_\nu)(1)\ ;$$

hence $M(s)f^\phi_{\tau,s}(1)$ has a pole (simple) at $s = 1$. Note that $f^\phi_{\tau,s}$ is a smooth section. Let $S_\infty$ be the set of archimedean places. Then

$$(\phi_\nu)_{\nu\in S_\infty} \mapsto \ell\Big(\otimes_{\nu\in S_\infty} M_\nu(1)\varphi_{\phi_\nu,1}(1)\otimes(\otimes_{\nu\notin S}\xi_\nu))(1)$$

is a continuous nontrivial linear functional on $\bigotimes_{\nu\in S_\infty} \mathrm{Ind}^{G_\nu}_{P_\nu} \tau_\nu$, and hence it is nontrivial on the dense subspace of $\prod_{\nu\in S_\infty} K_\nu$-finite vectors. Note, again, that the intertwining operator is holomorphic at $s = 1$ at each place (otherwise, the Eisenstein series has a high-order pole at $s = 1$, which is impossible.) This provides a pole at $s = 1$ for $M(s) f^{\phi}_{\tau,s}(1)$ with $\phi$ being $K$-finite. □

*Remark* 1. In the last proposition the choice of $S$ is immaterial. If there is a set $S_0$ such that $L^{S_0}(\tau, st, \frac{1}{2}) \neq 0$, then $L^S(\tau, st, \frac{1}{2}) \neq 0$ for any set $S$ (as above), and similarly, if $L^{S_0}(\tau, \Lambda^2, s)$ has a pole at $s = 1$, then $L^S(\tau, \Lambda^2, s)$ has a pole at $s = 1$, for any set $S$ (as above). The reason is that locally, at a finite place $\nu$, where $\tau_\nu$ is unramified $L(\tau_\nu, st, \frac{1}{2}) \neq 0$ and $L(\tau_\nu, \Lambda^2, 1) \neq 0$. Indeed, since $\tau_\nu$ is unitary the eigenvalues of its corresponding semisimple conjugacy class (in $\mathrm{GL}_{2n}(\mathbb{C})$) are all strictly less, in absolute value, then $q_\nu^{1/2}$ [J-S3] and so $L(\tau_\nu, st, s)$ (resp. $L(\tau_\nu, \Lambda^2, s)$) is holomorphic and nonzero at $s = \frac{1}{2}$ (resp. at $s = 1$). Similar reasoning in the archimedean case implies that our assumption on the standard $L$-function of $\tau$ is equivalent to $L(\tau, st, \frac{1}{2}) \neq 0$ (full $L$-function).

*Remark* 2. Let $\tau$ be a self-dual, irreducible, automorphic, cuspidal representation of either $\mathrm{GL}_{2n}(\mathbb{A})$ or $\mathrm{GL}_{2n+1}(\mathbb{A})$. As before, we can construct an Eisenstein series, induced from the Siegel-type parabolic subgroup of $G$, where $G$ is one of the following (split) groups

$$\mathrm{Sp}_{4n}, \mathrm{SO}_{4n}, \mathrm{SO}_{4n+1}, \widetilde{\mathrm{Sp}}_{4n+2}$$

(the last group is the metaplectic cover of $\mathrm{Sp}_{4n+2}$). The analogs of the quotient of (products of) $L$-functions in (1.3) are summarized in the following table (which includes the previous case of $\mathrm{Sp}_{4n}$)

(1.5)

| | |
|---|---|
| $\mathrm{Sp}_{4n}$ | $\frac{L^S(\tau,st,s-1/2)L^S(\tau,\Lambda^2,2s-1)}{L^S(\tau,st,s+\frac{1}{2})L^S(\tau,\Lambda^2,2s)}$ |
| $\mathrm{SO}_{4n}$ | $\frac{L^S(\tau,\Lambda^2,2s-1)}{L^S(\tau,\Lambda^2,2s)}$ |
| $\mathrm{SO}_{4n+1}$ | $\frac{L^S(\tau,Sym^2,2s-1)}{L^S(\tau,Sym^2,2s)}$ |
| $\widetilde{\mathrm{Sp}}_{4n+2}$ | $\frac{L^S(\tau,Sym^2,2s-1)}{L^S(\tau,Sym^2,2s)}$ |

.

Note that, in the last case, $\tau$ is a representation of $\mathrm{GL}_{2n+1}(\mathbb{A})$, and in the remaining cases, it is a representation of $\mathrm{GL}_{2n}(\mathbb{A})$. In all cases, except the last one, the functions $\phi$ and $f^{\phi}_{\tau,s}$ are defined in exactly the same way. In case $\widetilde{\mathrm{Sp}}_{4n+2}$, we require that $\phi$ is a smooth function on $\widetilde{\mathrm{Sp}}_{4n+2}(\mathbb{A})$, taking values in the space of $\tau$, such that

$$\phi((\widehat{m}u, 1)g; r) = \gamma_\psi(\det m)(\det m, \det m)\delta_P^{1/2}(\widehat{m})\phi(g; rm) \ , \tag{1.6}$$

for $m \in \mathrm{GL}_{2n+1}(\mathbb{A})$, $u \in U_{\mathbb{A}}$, $g \in \widetilde{\mathrm{Sp}}_{4n+2}(\mathbb{A})$, $r \in \mathrm{GL}_{2n+1}(\mathbb{A})$. Here $\gamma_\psi(x)$ is the (global) Weil factor associated to a nontrivial additive character $\psi$ of $F\backslash\mathbb{A}$. $\gamma_\psi$ is a character of the two-fold cover of $\mathbb{A}^*$. It satisfies $\gamma_\psi(x_1 \cdot x_2) = \gamma_\psi(x_1)\gamma_\psi(x_2)(x_1, x_2)$, for $x_1, x_2$ in $\mathbb{A}^*$. $(\,,\,)$ is the (global) Hilbert symbol. $f^{\phi}_{\tau,s}(g)$ is constructed by multiplying $\phi(g)$ by $H(g)^{s-1/2}$. We can now repeat the proof of the last proposition (word-for-word) and conclude that the corresponding Eisenstein series has a pole (simple) at $s = 1$, if $L^S(\tau, \Lambda^2, s)$ has a pole at $s = 1$ in case $\mathrm{SO}_{4n}$, or $L^S(\tau, \mathrm{Sym}^2, s)$ has a pole at $s = 1$ in cases $\mathrm{SO}_{4n+1}$ or $\widetilde{\mathrm{Sp}}_{4n+2}$.

2. *The* $\mathrm{Sp}_{2n} \times \mathrm{Sp}_{2n}$*-period of* $\mathrm{Res}_{s=1} E(g, f^{\phi}_{\tau,s})$. We go back to the case $G = \mathrm{Sp}_{4n}$, $\tau$ – self-dual, irreducible, automorphic, cuspidal representation. We assume that $L(\tau, st, \frac{1}{2}) \neq 0$ and that there exists a finite set of places $S$, including those at infinity, outside which $\tau$ is unramified, such that $L^S(\tau, \Lambda^2, s)$ has a pole at $s = 1$. (See Remark 1 in the last section.) Note that this implies, in particular, that $\tau$ has a trivial central character [J-S1]. By Proposition 1, $\mathrm{Res}_{s=1} E(g, f^{\phi}_{\tau,s})$ is nontrivial ($s = 1$ is a simple pole.) Consider the subgroup $H = \mathrm{Sp}_{2n} \times \mathrm{Sp}_{2n}$, embedded in $\mathrm{Sp}_{4n}$ by

$$\left( \begin{pmatrix} A_1 & B_1 \\ C_1 & D_1 \end{pmatrix}, \begin{pmatrix} A_2 & B_2 \\ C_2, & D_2 \end{pmatrix} \right) \overset{i}{\mapsto} \begin{pmatrix} A_1 & & & B_1 \\ & A_2 & -B_2 & \\ & -C_2 & D_2 & \\ C_1 & & & D_1 \end{pmatrix} . \tag{1.7}$$

Each letter represents an $n \times n$ block. We sometimes identify, to our convenience, $h$ and $i(h)$ in $\mathrm{Sp}_{4n}$, for $h \in H$. Denote

$$E_1(g, \phi) = \mathrm{Res}_{s=1} E(g, f^{\phi}_{\tau,s}) .$$

The main result of this chapter is

THEOREM 2. *Under the above assumption, $E_1$ is integrable over $H_F\backslash H_{\mathbb{A}}$ and (for a suitable choice of measures)*

$$\int_{H_F\backslash H_{\mathbb{A}}} E_1(h, \phi)dh = \int_{K_H} \int_{C_{2n}(\mathbb{A})\mathrm{GL}_n(F)^2\backslash \mathrm{GL}_n(\mathbb{A})^2} \phi(k; \begin{pmatrix} a & \\ & b \end{pmatrix})d(a, b)dk . \tag{1.8}$$

*Here $K_H = K \cap H$ and $C_{2n}$ is the center of* $\mathrm{GL}_{2n}$.

This theorem and its proof are very similar to Theorem 1 in [J-R1]. Note that the inner $\mathrm{GL}_n^2$-integration, on the right-hand side of (1.8), is the integral (32) in [F-J.] (with $s = \frac{1}{2}, \chi = \eta = 1$). By Theorem 4.1 in [F-J], a necessary condition, for this integral to be nonzero, is that $L^S(\tau, \Lambda^2, s)$ has a pole at $s = 1$, and, in this case, the integral is of the form $\alpha(k; \phi)L(\tau, st, \frac{1}{2})$, where $\alpha(k; \phi)$ is a nontrivial linear form, for each $k \in K_H$. Under our assumptions on $\tau$, the right-hand side of (1.8) is not identically zero. The proof for this is entirely similar to that of Proposition 2 in [J-R1]. All the requirements there are supplied by [F-J.]. Thus $\alpha(\phi) = \int_{K_H} \alpha(k,\phi)dk$ is nontrivial and hence the period of $E_1$ along $H$ is nontrivial. Note that, exactly as remarked on the bottom of p. 178 in [J-R1], formula (1.8) supplies another proof for the fact that if $L(\tau, st, \frac{1}{2}) \neq 0$ and $L^S(\tau, \Lambda^2, s)$ has a pole at $s = 1$, then $E(g, f^{\phi}_{\tau,s})$ has a pole at $s = 1$. Moreover, the last remarks, together with (1.8) prove

COROLLARY 3. *Let $\tau$ be a self-dual, irreducible, automorphic, cuspidal representation of $\mathrm{GL}_{2n}(\mathbb{A})$. Then $E(g, f^{\phi}_{\tau,s})$ has a pole at $s = 1$ and $E_1(\cdot, \phi)$ admits a nontrivial period along $H$, if and only if $L(\tau, st, \frac{1}{2}) \neq 0$ and $L^S(\tau, \Lambda^2, s)$ has a pole at $s = 1$.*

The rest of this chapter is devoted to the proof of Theorem 2.

3. *Truncation.* As in [J-R1], we consider the truncation operator applied to $E(g, f^{\phi}_{\tau,s})$. Denote, for a real number $c$, by $\chi_c$ the characteristic function of all real numbers larger than $c$. The only nontrivial constant terms of $E(g, f^{\phi}_{\tau,s})$ along unipotent radicals of standard parabolic subgroups, are those taken along $U$ or $\{1\}$. Thus, the truncation operator $\Lambda^c$, $c > 1$, applied to $E(g, f^{\phi}_{\tau,s})$ is

$$\Lambda^c E(g, f^{\phi}_{\tau,s}) = E(g, f^{\phi}_{\tau,s}) - \sum_{\gamma \in P_F \backslash G_F} E^U(\gamma g, f^{\phi}_{\tau,s}) \chi_c(H(\gamma g)) \ . \tag{1.9}$$

Since $E(g, f^{\phi}_{\tau,s})$ is smooth and of moderate growth, $\Lambda^c E(g, f^{\phi}_{\tau,s})$ is rapidly decreasing. Also, the sum on the right-hand side of (1.9) has finitely many terms (depending on $g$ and $c$). In particular, $\Lambda^c E(g, f^{\phi}_{\tau,s})$ is meromorphic in $s$. (See [A1], [A2], for more details.) Similar remarks are valid for $\Lambda^c E_1(g, \phi)$. By (1.1) and (1.9), we have

(1.10)

$$\Lambda^c E(g, f^{\phi}_{\tau,s}) = E(g, f^{\phi}_{\tau,s}) - \sum_{\gamma \in P_F \backslash G_F} (f^{\phi}_{\tau,s}(\gamma g) + M(s) f^{\phi}_{\tau,s}(\gamma g)) \chi_c(H(\gamma g))$$

$$= \sum_{\gamma \in P_F \backslash G_F} f^{\phi}_{\tau,s}(\gamma g) \chi^c(H(\gamma g)) - \sum_{\gamma \in P_F \backslash G_F} M(s) f^{\phi}_{\tau,s}(\gamma g) \chi_c(H(\gamma g)) \ .$$

This last equality is for $\mathrm{Re}(s) > n + 1$. $\chi^c$ is the characteristic function of the interval $(0, c]$. Denote the first sum, in the last expression of (1.10), by $\theta^c_1(g, f^{\phi}_{\tau,s})$, and the second – by $\theta^c_2(g, f^{\phi}_{\tau,s})$, so that

$$\Lambda^c E(g, f^{\phi}_{\tau,s}) = \theta^c_1(g, f^{\phi}_{\tau,s}) - \theta^c_2(g, f^{\phi}_{\tau,s}) \ . \tag{1.11}$$

Note again that the sum defining $\theta^c_2(g, f^{\phi}_{\tau,s})$ has finitely many terms (depending on $g$ and on $c$). In particular $\theta^c_2(g, f^{\phi}_{\tau,s})$ is meromorphic in $s$, and so $\theta^c_1(g, f^{\phi}_{\tau,s})$ defines a meromorphic function in the whole plane. (1.11) is now valid as an equality of meromorphic functions.

Denote $\mathrm{Res}_{s=1} M(s) = M_1$. Since $f^{\phi}_{\tau,s}$ is holomorphic, then an application of $\Lambda^c$ to $E_1$ yields similarly,

$$\Lambda^c E_1(g,\phi) = E_1(g,\phi) - \theta^c_3(g,\phi) \ , \tag{1.12}$$

where

$$\theta^c_3(g,\phi) = \sum_{\gamma\in P_F\backslash G_F} M_1(f^{\phi}_{\tau,1}(\gamma g)\chi_c(H(\gamma g)) \ . \tag{1.13}$$

Since $\Lambda^c E(g, f^{\phi}_{\tau,s})$ and $\Lambda^c E_1(g,\phi)$ are rapidly decreasing, they are bounded and hence integrable on $H_F\backslash H_{\mathbb{A}}$. We will prove

Proposition 4. *$\theta^c_1(\cdot, f^{\phi}_{\tau,s})$ is integrable on $H_F\backslash H_{\mathbb{A}}$, if $\mathrm{Re}(s)$ is sufficiently large, $\theta^c_2(\cdot, f^{\phi}_{\tau,s})$ is integrable on $H_F\backslash H_{\mathbb{A}}$, if $\mathrm{Re}(s) > 0$, (and $M(s)$ exists) and $\theta^c_3(\cdot,\phi)$ is integrable on $H_F\backslash H_{\mathbb{A}}$. The following formulae are valid, with a certain choice of measures*:

(1.14)
$$\int\limits_{H_F\backslash H_{\mathbb{A}}} \theta^c_1(h, f^{\phi}_{\tau,s})dh = \frac{c^{s-1}}{s-1}\int\limits_{K_H}\ \int\limits_{C_{2n}(\mathbb{A})\mathrm{GL}_n(F)^2\backslash \mathrm{GL}_n(\mathbb{A})^2} \phi(k;\begin{pmatrix} a & \\ & b\end{pmatrix}d(a,b)dk \ ,$$

(1.15)
$$\int\limits_{H_F\backslash H_{\mathbb{A}}} \theta^c_2(h, f^{\phi}_{\tau,s})dh$$
$$= \frac{c^{-s}}{s}\int\limits_{K_H}\ \int\limits_{C_{2n}(\mathbb{A})\mathrm{GL}_n(F)^2\backslash \mathrm{GL}_n(\mathbb{A})^2} M(s)(\varphi^{\phi}_{\tau,s})(k;\begin{pmatrix} a & \\ & b\end{pmatrix})d(a,b)dk \ ,$$

(1.16)
$$\int\limits_{H_F\backslash H_{\mathbb{A}}} \theta^c_3(h,\phi)dh$$
$$= c^{-1}\int\limits_{K_H}\ \int\limits_{C_{2n}(\mathbb{A})\mathrm{GL}_n(F)^2\backslash \mathrm{GL}_n(\mathbb{A})^2} M_1(\varphi^{\phi}_{\tau,1}(k;\begin{pmatrix} a & \\ & b\end{pmatrix}d(a,b)dk \ .$$

This proposition will finish the proof of Theorem 2, exactly as in [J-R1, p. 181]. Integrating (1.11) along $H_F\backslash H_{\mathbb{A}}$, first for $\mathrm{Re}(s) \gg 0$, and using (1.14) and (1.15), we get

$$\int\limits_{H_F\backslash H_{\mathbb{A}}} \Lambda^c E(h, f^{\phi}_{\tau,s})dh = \frac{c^{s-1}}{s-1}\iint \phi(k;\begin{pmatrix} a & \\ & b\end{pmatrix}d(a,b)dk \tag{1.17}$$
$$- \frac{c^{-s}}{s}\iint M(s)(\varphi^{\phi}_{\tau,s})(k;\begin{pmatrix} a & \\ & b\end{pmatrix})d(a,b)dk \ .$$

By analytic continuation, this makes sense in the whole plane.

Taking residues in (1.17) at $s = 1$,

(1.18)
$$\begin{aligned}\int\limits_{H_F\backslash H_{\mathbb{A}}} \Lambda^c E_1(h,\phi)dh &= \iint \phi(k;\begin{pmatrix} a & \\ & b\end{pmatrix}) d(a,b)dk \\ &\quad - c^{-1}\iint M_1(\varphi^{\phi}_{\tau,1})(k;\begin{pmatrix} a & \\ & b\end{pmatrix})d(a,b)dk \\ &= \iint \phi(k;\begin{pmatrix} a & \\ & b\end{pmatrix})d(a,b)dk - \int\limits_{H_F\backslash H_{\mathbb{A}}} \theta^c_3(h,\phi)dh\ .\end{aligned}$$

We used (1.16). Comparing (1.18) with (1.12), we conclude that $E_1$ is integrable on $H_F\backslash H_{\mathbb{A}}$ and that (1.8) is satisfied.

Each of the functions $\theta^c_j$ has the form

$$\theta_j(g) = \sum_{\gamma\in P_F\backslash G_F} \xi_j(\gamma g)\ . \tag{1.19}$$

where

$$\xi_j(g) = \begin{cases} f^{\phi}_{\tau,s}(g)\chi^c(H(g))\ , & j=1 \\ M(s)(f^{\phi}_{\tau,s})(g)\chi_c(H(g))\ , & j=2 \\ M_1(f^{\phi}_{\tau,1})(g)\chi_c(H(g))\ , & j=3\end{cases}\ . \tag{1.20}$$

The sum (1.19), in case $j = 1$, converges only for $\mathrm{Re}(s) \gg 0$. Before proving the integrability of $\theta_j$ on $H_F\backslash H_{\mathbb{A}}$, and then compute its period on $H$, we proceed formally,

(1.21)
$$\begin{aligned}\int\limits_{H_F\backslash H_{\mathbb{A}}} \theta_j(h)dh &= \int\limits_{H_F\backslash H_{\mathbb{A}}} \sum_{\gamma\in P_F\backslash G_F} \xi_j(\gamma h)dh \\ &= \sum_{\gamma\in P_F\backslash G_F/H_F} \int\limits_{\gamma^{-1}P_F\gamma\cap H_F\backslash H_{\mathbb{A}}} \xi_j(\gamma h)dh\ .\end{aligned}$$

The set $P_F\backslash G_F/H_F$ is finite and will be described soon. Our task will be to show the integrability of each $\xi_j(\gamma h)$ (in $h$)) on $\gamma^{-1}P_F\gamma \cap H_F\backslash H_{\mathbb{A}}$, and then compute the integral.

4. *The set $P_F\backslash G_F/H_F$.* The description of this set is known. It appears as an ingredient in the "doubling method." See [PS-R, Lemma 2.1], for the description of the action of $H_F$ on $P_F\backslash G_F$, realized as the variety of maximal isotropic subspaces of the symplectic space of dimension $4n$. Using this description and passing to a different basis of $F^{4n}$, so that $H = \mathrm{Sp}_{2n} \times \mathrm{Sp}_{2n}$ is embedded in $G$ by (1.7), we arrive at the following set of representatives for $P_F\backslash G_F/H_F$.

$$\gamma_d = \begin{pmatrix} I_d & & \\ & \gamma'_d & \\ & & I_d \end{pmatrix} , \quad 0 \le d \le n , \tag{1.22}$$

where

$$\gamma'_d = \left(\begin{array}{ccc|ccc} & & I_{n-d} & & & \\ & I_d & & & & \\ 0 & & & & & -I_{n-d} \\ \hline I_{n-d} & & -I_{n-d} & & & 0 \\ & & & & I_d & \\ & & & I_{n-d} & & I_{n-d} \end{array}\right) . \tag{1.23}$$

Note that

$$\gamma'_d = \left(\begin{array}{ccc|ccc} & & I_{n-d} & & & \\ & I_d & & & & \\ 0 & & & & & -I_{n-d} \\ \hline I_{n-d} & & & & & 0 \\ & & & & I_d & \\ & & & I_{n-d} & & \end{array}\right) \left(\begin{array}{ccc|ccc} I_{n-d} & & -I_{n-d} & & & \\ & I_d & & & & \\ & & I_{n-d} & & & \\ \hline & & & I_{n-d} & & I_{n-d} \\ & & & & I_d & \\ & & & & & I_{n-d} \end{array}\right) .$$

A computation of the stabilizer of $\gamma_d$ shows that

(1.24)

$$H_F \cap \gamma_d^{-1} P_F \gamma_d$$
$$= \left\{ \left( \begin{pmatrix} a & x & y \\ & c & x' \\ & & a^* \end{pmatrix} , \begin{pmatrix} b & u & z \\ & c & u' \\ & & b^* \end{pmatrix} \right) \in H_F \,\middle|\, \begin{matrix} a, b \in \mathrm{GL}_d(F) \\ c \in \mathrm{Sp}_{2(n-d)}(F) \end{matrix} \right\} .$$

Write this as the semidirect product $M_d \ltimes V_d$, where

$$M_d = \left\{ \left( \begin{pmatrix} a & & \\ & c & \\ & & a^* \end{pmatrix} , \begin{pmatrix} b & & \\ & c & \\ & & b^* \end{pmatrix} \right) \,\middle|\, \begin{matrix} a, b \in \mathrm{GL}_d(F) \\ c \in \mathrm{Sp}_{2(n-d)}(F) \end{matrix} \right\}$$

$$V_d = \left\{ \left( \begin{pmatrix} I_d & x & y \\ & I_{2(n-d)} & x' \\ & & I_d \end{pmatrix} , \begin{pmatrix} I_d & u & z \\ & I_{2(n-d)} & u' \\ & & I_d \end{pmatrix} \right) \in H_F \right\} .$$

We now collect several conjugation formulae which will be used later. For

$$r = i\left(\begin{pmatrix} a & & \\ & I_{2(n-d)} & \\ & & a^* \end{pmatrix}, \begin{pmatrix} b & & \\ & I_{2(n-d)} & \\ & & b^* \end{pmatrix}\right),$$

($i$ is the embedding (1.7))

$$\gamma_d r \gamma_d^{-1} = \begin{pmatrix} a &&&&&&& \\ & I_{n-d} &&&&&& \\ && b &&&&& \\ &&& I_{n-d} &&&& \\ &&&& I_{n-d} &&& \\ &&&&& b^* && \\ &&&&&& I_{n-d} & \\ &&&&&&& a^* \end{pmatrix}. \tag{1.25}$$

For $r = i\left(\begin{pmatrix} I_d & & \\ & c & \\ & & I_d \end{pmatrix}, \begin{pmatrix} I_d & & \\ & c & \\ & & I_d \end{pmatrix}\right)$, where $c = \begin{pmatrix} c_1 & c_2 \\ c_3 & c_4 \end{pmatrix}$

is the $(n-d) \times (n-d)$ block description of $c \in \mathrm{Sp}_{2(n-d)}(F)$,

$$\gamma_d r \gamma_d^{-1} = \left(\begin{array}{cccc|cccc} I_d &&&&&&& 0 \\ & c_1 && -c_2 &&& -c_2 & \\ && I_d &&& 0 && \\ & -c_3 && c_4 & -c_3 &&& \\ \hline &&&& c_1 && c_2 & \\ &&&&& I_d && \\ &&&& c_3 && c_4 & \\ &&&&&&& I_d \end{array}\right). \tag{1.26}$$

For

$$v = i\left(\begin{pmatrix} I_d & x_1 & e_1 & y_1 \\ & I_{n-d} & 0 & e_1' \\ & & I_{n-d} & -x_1' \\ & & & I_d \end{pmatrix}, \begin{pmatrix} I_d & x_2 & e_2 & y_2 \\ & I_{n-d} & 0 & e_2' \\ & & I_{n-d} & -x_1' \\ & & & I_d \end{pmatrix}\right) \in V_d,$$

$$(1.27) \qquad \gamma_d v \gamma_d^{-1} = \left(\begin{array}{cccc|cccc} I_d & x_1 & 0 & -e_1 & x_1 & 0 & 0 & y_1 \\ & I_{n-d} & 0 & 0 & & -e_2' & 0 & 0 \\ & x_2 & I_d & -e_2 & & -y_2 & -e_2 & 0 \\ & & & I_{n-d} & & & & x_1' \\ \hline & & & & I_{n-d} & e_2' & 0 & e_1' \\ & & & & & I_d & 0 & 0 \\ & & & & & -x_2' & I_{n-d} & -x_1' \\ & & & & & & & I_d \end{array}\right) .$$

For $3m' = i\left(\begin{pmatrix} I_d & & & \\ & e & ez & \\ & 0 & e^* & \\ & & & I_d \end{pmatrix}, I_{2n}\right)$,

$$(1.28) \qquad \gamma_d m' \gamma_d^{-1} = \left(\begin{array}{cccc|cccc} I_d & & & & & & & \\ & I_{n-d} & & & & & & \\ & & I_d & & & & & \\ & & & e^* & & & & \\ \hline & -I_{n-d} & 0 & -ez & e & & & \\ & & & 0 & & I_d & & \\ & & & -e^* & & & I_{n-d} & \\ & & & 0 & & & & I_d \end{array}\right) k_0$$

where

$$k_0 = \left(\begin{array}{cccc|cccc} I_d & & & & & & & \\ & I_{n-d} & & & & & & \\ & & I_d & & & & & \\ & & & I_{n-d} & & & & \\ \hline & I_{n-d} & 0 & 0 & I_{n-d} & & & \\ & & & 0 & & I_d & & \\ & & & I_{n-d} & & & I_{n-d} & \\ & & & 0 & & & & I_d \end{array}\right) .$$

5. *A formal proof of Proposition* 4. We first prove Proposition 4 formally, without paying attention to convergence issues. Denote

$$Q_d = \gamma_d^{-1} P \gamma_d \cap H \ ,$$

$$Q_d' = \left\{ \left( \begin{pmatrix} a & x & y \\ & c & x' \\ & & a^* \end{pmatrix} , \begin{pmatrix} b & u & z \\ & e & u' \\ & & b^* \end{pmatrix} \right) \in H \,\middle|\, \begin{matrix} a, b \in \mathrm{GL}_d \\ c, e \in \mathrm{Sp}_{2(n-d)} \end{matrix} \right\} .$$

$Q'_d$ is a parabolic subgroup. Write its Levi decomposition as $M'_d \ltimes V_d$. Note that $M_d \subset M'_d$.
Following (1.21), we have to compute the integrals

$$I_{j,d} = \int\limits_{Q_d(F)\backslash H_{\mathbb{A}}} \xi_j(\gamma_d h) \quad \text{for} \quad 0 \le d \le n \quad \text{and} \quad j = 1,2,3 \ .$$

Let $0 < d < n$. Write the Iwasawa decomposition in $H_{\mathbb{A}}$,

$$h = vm'k \ , \qquad m' \in M'_d(\mathbb{A}) \ , \qquad v \in V_d(\mathbb{A}) \ , \qquad k \in K_H \ ,$$
$$dh = \delta^{-1}(m')dvdm'dk \ ,$$

where $\delta$ is the modular function corresponding to the parabolic subgroup $Q'_d$. Then

(1.29)
$$I_{j,d} = \int\limits_{K_H} \int\limits_{M_d(F)\backslash M'_d(\mathbb{A})} \int\limits_{V_d(F)\backslash V_d(\mathbb{A})} \xi_j(\gamma_d v \gamma_d^{-1} \cdot \gamma_d m'k)\delta^{-1}(m')dvdm'dk \ .$$

By (1.27), the projection to $\mathrm{GL}_{2n}$ of $\gamma_d v \gamma_d^{-1}$, as $v$ varies in $V_d$ is a unipotent radical $N_d$ in $\mathrm{GL}_{2n}$, (since $0 < d < n$). Then $H(\gamma_d v \gamma_d^{-1}) = 1$, and from (1.20) we see that (1.29) involves an inner integration of a cusp form in $\tau$ along $N_d(F)\backslash N_d(\mathbb{A})$, and hence

$$(1.30) \qquad I_{j,d} = 0 \ , \qquad \text{for} \quad 0 < d < n, \ j = 1,2,3 \ .$$

In case $d = 0$, $Q_0 = \{(c,c) \mid c \in \mathrm{Sp}_{2n}\}$. We have

$$(1.31) \qquad I_{j,0} = \int\limits_{\mathrm{Sp}_{2n}(\mathbb{A})} \int\limits_{\mathrm{Sp}_{2n}(F)\backslash \mathrm{Sp}_{2n}(\mathbb{A})} \xi_j\Big(\gamma_0 i(c,c)\gamma_0^{-1}\gamma_0 i(e,1)\Big) dcde \ .$$

By (1.26), the projection of $\gamma_0 i(b,b)\gamma_0^{-1}$ on $\mathrm{GL}_{2n}$, as $b$ varies in $\mathrm{Sp}_{2n}$ is $\mathrm{Sp}_{2n}$ (embedded as a subgroup of $\mathrm{GL}_{2n}$). Since $H(\gamma_0 i(c,c)\gamma_0^{-1}) = 1$, we obtain in (1.31) an inner integration of a cusp form in $\tau$ along $\mathrm{Sp}_{2n}(F)\backslash \mathrm{Sp}_{2n}(\mathbb{A})$, which equals zero by [J-R1, Prop. 1]. Thus,

$$(1.32) \qquad I_{j,0} = 0 \ , \qquad \text{for} \quad j = 1,2,3 \ .$$

The only remaining case is $d = n$. Here $\gamma_n = I_{4n}$ and

$$Q_n = \left\{ \begin{pmatrix} a & y \\ 0 & a^* \end{pmatrix}, \begin{pmatrix} b & z \\ 0 & b^* \end{pmatrix} \in H \mid a, b, \in \mathrm{GL}_n \right\} = Q'_n \ .$$

Using (1.25) and (1.27), we have

(1.33)

$$I_{j,n} = \int\limits_{K_H} \int\limits_{\mathrm{GL}_n(F)^2\backslash \mathrm{GL}_n(\mathbb{A})^2} \xi_j\left(\begin{pmatrix} a &&& \\ & b && \\ && b^* & \\ &&& a^* \end{pmatrix} k\right) |\det(ab)|^{-(n+1)} d(a,b)dk$$

$$= \int\limits_{K_H} \int\limits_{\mathrm{GL}_n(F)^2\backslash (\mathrm{GL}_n(\mathbb{A})^2)^0} \int\limits_{F^*\backslash \mathbb{A}^*} \xi_j\left(\begin{pmatrix} ta &&& \\ & tb && \\ && t^{-1}b^* & \\ &&& t^{-1}a^* \end{pmatrix} k\right)$$

$$\cdot\ |t|^{-2n(n+1)} |\det(ab)| d^*t d(a,b)dk\ .$$

Here $(\mathrm{GL}_n(A^2))^0 = \{(a,b) \in \mathrm{GL}_n(\mathbb{A})^2 |\ |\det(ab)| = 1\}$. Of course, in (1.33), we made a choice of measures. $d^*t$ will soon be specified. Thus, by (1.20) (recall also that, by our assumptions, $\tau$ has a trivial central character)

(1.34)

$$I_{1,n} = \int\limits_{K_H} \int\limits_{\mathrm{GL}_n(F)\backslash (\mathrm{GL}_n(\mathbb{A})^2)^0} \phi(k; \begin{pmatrix} a & \\ & b \end{pmatrix}) d(a,b)dk \int\limits_{\substack{t\in F^*\backslash \mathbb{A}^* \\ |t|^{2n}\le c}} |t|^{2n(s-1)} d^*t\ .$$

Choose $d^*t$ such that the corresponding $d^*t$-integral equals $\int\limits_0^c t^{s-1}\frac{dt}{t} = \frac{c^{s-1}}{s-1}$, for $\mathrm{Re}(s) > 1$. The $d(a,b)$-integration yields one on $C_{2n}(\mathbb{A})\mathrm{GL}_n(F)^2\backslash \mathrm{GL}_n(\mathbb{A})^2$. This proves (1.14). In a similar fashion, we get, for $\mathrm{Re}(s) > 0$,

(1.35)

$$I_{2,n} = \int\limits_{K_H} \int\limits_{\mathrm{GL}_n(F)^2\backslash (\mathrm{GL}_n(\mathbb{A})^2)^0} M(s)(\varphi^{\phi}_{\tau,s})(k; \begin{pmatrix} a & \\ & b \end{pmatrix}) d(a,b)dk \int\limits_c^\infty t^{(1-s)-1}\frac{dt}{t}$$

$$= \frac{c^{-s}}{s} \int\limits_{K_H} \int\limits_{C_{2n}(\mathbb{A})\mathrm{GL}_n(F)^2\backslash \mathrm{GL}_n(\mathbb{A})^2} M(s)(\varphi^{\phi}_{\tau,s})(k; \begin{pmatrix} a & \\ & b \end{pmatrix}) d(a,b)dk\ ,$$

and

(1.36)

$$I_{3,n} = c^{-1} \int\limits_{K_H} \int\limits_{C_{2n}(\mathbb{A})\mathrm{GL}_n(F)^2\backslash \mathrm{GL}_n(\mathbb{A})^2} M_1(\varphi^{\phi}_{\tau,1})(k; \begin{pmatrix} a & \\ & b \end{pmatrix}) d(a,b)dk\ .$$

In the next section, we justify our formal calculations.

6. *Justification of the formal proof of Proposition* 4. By (1.21), we have to examine the absolute convergence of the integrals $I_{j,d}$. The case $d = n$ is, actually, rigorously proved in the end of the last section (1.34)–(1.36): $I_{1,n}$ (resp. $I_{2,n}$) converges absolutely for $\mathrm{Re}(s) > 1$ (resp. $\mathrm{Re}(s) > 0$) and $I_{3,n}$ converges absolutely.

PROPOSITION 5. *Let* $0 \le d < n$. *Then*:

1. $I_{1,d}$ *converges absolutely, for* $\mathrm{Re}(s)$ *large enough.*
2. $I_{2,d}$ *converges absolutely, for all* $s$ *where* $M(s)$ *is defined.*
3. $I_{3,d}$ *converges absolutely.*

*In all three cases,* $I_{j,d} = 0$.

*Proof.* We start with formula (1.29). We will prove the absolute convergence (at the indicated domains in $s$) of

(1.37)
$$\int\limits_{K_H} \int\limits_{M_d(\mathbb{A})\backslash M'_d(\mathbb{A})} \int\limits_{M_d(F)\backslash M_d(\mathbb{A})} \int\limits_{V_d(F)\backslash V_d(\mathbb{A})} \cdot\ \xi_j(\gamma_d v \gamma_d^{-1} \cdot \gamma_d r \gamma_d^{-1} \cdot \gamma_d m' k)\delta^{-1}(r)dvdrdm'dk.$$

We realize $M_d(\mathbb{A})\backslash M'_d(\mathbb{A})$ as $\left\{ \left( \begin{pmatrix} I_d & & \\ & x & \\ & & I_d \end{pmatrix} , I_{2n} \right) \Big| x \in \mathrm{Sp}_{2(n-d)}(\mathbb{A}) \right\}$ and then use the Iwasawa decomposition, following the Siegel parabolic subgroup $R' = L' \ltimes U'$, in $\mathrm{Sp}_{2(n-d)}$.

$$x = \begin{pmatrix} e & ez \\ & e^* \end{pmatrix} k'$$
$$dx = dedzdk' \ .$$

Here $e \in \mathrm{GL}_{n-d}(\mathbb{A})$, $z \in \mathcal{U}'_{\mathbb{A}} = \left\{ z | \begin{pmatrix} I_{n-d} & z \\ & I_{n-d} \end{pmatrix} \in U'_{\mathbb{A}} \right\}$, $k' \in K'$ — the standard maximal compact subgroup of $\mathrm{Sp}_{2(n-d)}(\mathbb{A})$. Let

$$\omega_d = \begin{pmatrix} I_d & & & \\ & 0 & I_d & \\ & I_{n-d} & 0 & \\ & & & I_{n-d} \end{pmatrix} .$$

Then, by (1.25)–(1.28), the integral (1.37) equals

$$\int \phi_j(k_{e,z}k_0\gamma_d k'k;\begin{pmatrix} I_{2d} & y \\ & I_{2(n-d)}\end{pmatrix}\begin{pmatrix} a & & \\ & b & \\ & & c\end{pmatrix}\omega_d r_{e,z})|ab|^{s_j}|r_{e,z}|^{s'_j} \cdot \chi_{j,c}(|ab|\ |r_{e,z}|)dyd(a,b)dcdedzdk'dk\ . \tag{1.38}$$

The integration is on

$$y\in M_{2d\times 2(n-d)}(F)\backslash M_{2d\times 2(n-d)}(\mathbb{A}),\ \ (a,b)\in \mathrm{GL}_d(F)^2\backslash \mathrm{GL}_d(\mathbb{A})^2\ ,$$
$$c\in \mathrm{Sp}_{2(n-d)}(F)\backslash \mathrm{Sp}_{2(n-d)}(\mathbb{A}),\ \ e\in \mathrm{GL}_{n-d}(\mathbb{A}),\ \ z\in \mathcal{U}'_{\mathbb{A}},\ \ k'\in K'$$

(embedded naturally in $\mathrm{Sp}_{4n}(\mathbb{A})$)and $k\in K_H$. $k_{e,z}$ is the compact part of the Iwasawa decomposition in $\mathrm{Sp}_{4n}(\mathbb{A})$ of the matrix, multiplying $k_0$, on the right-hand side of (1.28), and $r_{e,z}$ is the projection to $\mathrm{GL}_{2(n-d)}(\mathbb{A})$ of its "Siegel parabolic part." $\phi_j$ is defined as follows (see (1.20)). For $g\in G_{\mathbb{A}}$ and $m\in \mathrm{GL}_{2n}(\mathbb{A})$,

$$\phi_j(g;m)=\begin{cases}\varphi^{\phi}_{\tau,s}(g;m)\ , & j=1\\ M(s)(\varphi^{\phi}_{\tau,s})(g;m)\ , & j=2\\ M_1(\varphi^{\phi}_{\tau,1})(g;m)\ , & j=3\end{cases}\ ,$$

and for $t\in\mathbb{R}$

$$\chi_{j,c}(t)=\begin{cases}\chi^c(t) & j=1\\ \chi_c(t) & j=2,3\end{cases}\ .$$

Finally,

$$s_j=\begin{cases} s-n+d-1\ , & j=1\\ 1-s-n+d-1\ , & j=2\\ -n+d-1\ , & j=3\end{cases}\ ,$$

$$s'_j=\begin{cases} s+n\ , & j=1\\ 1-s+n\ , & j=2\\ n\ , & j=3\end{cases}\ .$$

In (1.38), we denote, for a matrix $x$, $|x|=|\det x|$. Examining (1.28), we see that

$$\omega_d r_{e,z}=m'_{e,z}\omega_d\ , \tag{1.39}$$

where

$$m'_{e,z}=\begin{pmatrix} I_{n+d} & \\ & e^*\end{pmatrix}\begin{pmatrix} I_{2d} & \\ & m_{e,z}\end{pmatrix}\ , \tag{1.40}$$

and $m_{e,z}$ is the projection to $\mathrm{GL}_{2(n-d)}(\mathbb{A})$ of the "Siegel parabolic part" in the Iwasawa decomposition in $\mathrm{Sp}_{2(n-d)}(\mathbb{A})$ of

$$u_{e,z} = \left(\begin{array}{cc|cc} I_{n-d} & & & \\ & I_{n-d} & & \\ \hline -e^{-1} & -z & I_{n-d} & \\ 0 & -e^* & & I_{n-d} \end{array}\right) .$$

We want to establish the absolute convergence of (1.38). For this, we use the majorization of cusp forms on $\mathrm{GL}_{2n}(\mathbb{A})$, given in Proposition 6 of [J-R1]. Since $\phi_j$ is $K$-finite, we find that, for each integer $N > 1$, there is a constant $c'$, such that

$$(1.41)\qquad \left|\phi_j(\widetilde{k}, \begin{pmatrix} m & y \\ 0 & r \end{pmatrix})\right| \leq c' \max\left\{1, |m|^{2(n-d)}|r|^{-2d}\right\}^{-N} ,$$

for all $\widetilde{k} \in K$, $m \in \mathrm{GL}_{2d}(\mathbb{A})$, $r \in \mathrm{GL}_{2(n-d)}(\mathbb{A})$. Now, it suffices to replace the $d(a,b)$-integration in (1.38) by an integration on $\mathrm{GL}_{2d}(F)\backslash\mathrm{GL}_{2d}(\mathbb{A})$. By (1.40) and (1.41), it remains to consider

$$(1.42)\qquad I'_{j,d} = \int \max\{1, |g|^{2(n-d)}(|e^*|\ |m_{e,z}|)^{-2d}\}^{-N}|g|^{s_j}(|e^*|\ |m_{e,z}|)^{s'_j} \cdot \chi_{j,c}(|g|\ |e^*|\ |m_{e,z}|)dgdedz ,$$

where $g$ is integrated over $\mathrm{GL}_{2d}(F)\backslash\mathrm{GL}_{2d}(\mathbb{A})$, $e$-over $\mathrm{GL}_{n-d}(\mathbb{A})$ and $z$ – over $\mathcal{U}'_{\mathbb{A}}$. Before proceeding to each case, let us prove two lemmas.

Lemma 6. *For all* $e \in \mathrm{GL}_{n-d}(\mathbb{A})$, $z \in \mathcal{U}'_{\mathbb{A}}$,

$$(1.43)\qquad |e^*|\ |m_{e,z}| \leq 1 .$$

*Proof.* Write the Iwasawa decomposition of $u_{e,z}$, at each place $\nu$, following the Borel subgroup of $\mathrm{Sp}_{2(n-d)}(F_\nu)$,

$$\overline{u}_\nu = u_{e_\nu,z_\nu} = \begin{pmatrix} t_1 & & & & & * \\ & \ddots & & & & \\ & & t_{2(n-d)} & & & \\ & & & t_{2(n-d)}^{-1} & & \\ & & & & \ddots & \\ & & & & & t_1^{-1} \end{pmatrix} k_\nu .$$

Then

$$(1.44)\qquad |m_{e_\nu,z_\nu}|^{-1} = |t_{2(n-d)} \cdot\ \cdots\ \cdot t_1|^{-1} = \|\varepsilon_{-2(n-d)}\overline{u}_\nu \wedge \cdots \wedge \varepsilon_{-1}\overline{u}_\nu\| ,$$

where $\{\varepsilon_1, \ldots, \varepsilon_{2(n-d)}; \varepsilon_{-2(n-d)}, \ldots, \varepsilon_{-2}, \varepsilon_{-1}\}$ is the standard basis of row vectors in $F_\nu^{4(n-d)}$. For a row vector $v = (x_1, \ldots, x_\ell) \in F_\nu^\ell$, we let

$$\|v\| = \begin{cases} \max_{1\le i\le \ell} |x_i| \ , & \nu < \infty \\ \big(\sum\limits_{i=1}^{\ell} x_i^2\big)^{1/2} \ , & F_\nu = \mathbb{R} \\ \sum\limits_{i=1}^{\ell} |x_i|^2 \ , & F_\nu = \mathbb{C} \ . \end{cases}$$

(If $\nu$ is archimedean, the norm on $\wedge^{2(n-d)} F_\nu^{4(n-d)}$ is defined through the orthonormal basis obtained from the standard basis $\{\varepsilon_1, \ldots, \varepsilon_{-1}\}$ above.) From (1.44), we get (with self-explanatory notation)

$$|m_{e_\nu, z_\nu}|^{-1} \ge \|(\varepsilon_{-(n-d)} \oplus E_{n-d}) \wedge \cdots \wedge (\varepsilon_{-1} \oplus E_1)\| \ .$$

where $E_{n-d}, \ldots E_1$ are the rows of the matrix $-e_\nu^*$;

$$-e_\nu^* = \begin{pmatrix} E_{n-d} \\ \vdots \\ E_1 \end{pmatrix} \ .$$

(See definition of $u_{e,z}$.) In particular,

$$|m_{e_\nu, z_\nu}|^{-1} \ge \|E_{n-d} \wedge \cdots \wedge E_1\| = |e_\nu^*| = |e_\nu|^{-1} \ .$$

This implies (1.43) at all places $\nu$; hence (1.43) is satisfied. This proves the lemma. ☐

Lemma 7. *Let $t$ be a real number. Then the integral*

$$A = \int\limits_{\mathcal{U}'_{\mathbb{A}}} \int\limits_{\mathrm{GL}_{n-d}(\mathbb{A})} (|e^*| \ |m_{e,z}|)^t dedz$$

*converges, provided $t$ is large enough.*

*Proof.* We first consider the local analog of $A$, at each place $\nu$,

$$A_\nu = \int\limits_{\mathcal{U}'_\nu} \int\limits_{\mathrm{GL}_{n-d}(F_\nu)} (|e^*| \ |m_{e,z}|)^t dedz \ .$$

Split $A_\nu$ into two summands $A_\nu^{(1)} + A_\nu^{(2)}$, where in $A_\nu^{(1)}$, we integrate over $|e| \le 1$ and in $A_\nu^{(2)}$ – over $|e| > 1$. We have

$$A_\nu^{(1)} = \iint\limits_{|e|\le 1} (|e^*| \ |m_{e,z}))^t dedz = \iint\limits_{|e|\ge 1} (|e| \ |\widetilde{m}_{e,z}|)^t dedz \ .$$

We changed variable $e \mapsto -e^{-1}$ and $z \mapsto -z$ and $\widetilde{m}_{e,z}$ is the projection to $\mathrm{GL}_{n-d}(\mathbb{A})$ of the "Siegel parabolic part" in the Iwasawa decomposition in $\mathrm{Sp}_{2(n-d)}(\mathbb{A})$ of

$$\widetilde{u}_{e,z} = \left(\begin{array}{cc|cc} I_{n-d} & & & \\ & I_{n-d} & & \\ \hline e & z & I_{n-d} & \\ e^{-1} & e' & & I_{n-d} \end{array}\right) .$$

Let $d_a(e)$ be the additive measure on $M_{n-d}(F_\nu)$. Then $de = \frac{d_a(e)}{|e|^{n-d}}$ and

$$A_\nu^{(1)} \leq \iint\limits_{|e|\geq 1} (|e|\ |\widetilde{m}_{e,z}|)^t d_a(e)dz \ . \tag{1.45}$$

As in (1.44),

$$\begin{aligned} |\widetilde{m}_{e,z}|^{-1} = \|&(\varepsilon_{-2(n-d)} \oplus \mathcal{E}_{n-d} \oplus z_{n-d}) \wedge \cdots \wedge (\varepsilon_{-(n-d)-1} \oplus \mathcal{E}_1 \oplus z_1) \\ &\wedge (\varepsilon_{-(n-d)} \oplus E_{n-d}) \wedge \cdots \wedge (\varepsilon_{-1} \oplus E_1)\| \ , \end{aligned} \tag{1.46}$$

where

$$e = \begin{pmatrix} \mathcal{E}_{n-d} \\ \vdots \\ \mathcal{E}_1 \end{pmatrix} \ , \ z = \begin{pmatrix} z_{n-d} \\ \vdots \\ z_1 \end{pmatrix} \ , \ e' = \begin{pmatrix} E_{n-d} \\ \vdots \\ E_1 \end{pmatrix} \ .$$

Hence

$$\begin{aligned} |\widetilde{m}_{e,z}|^{-1} &\geq \|(\varepsilon_{-2(n-d)} + \mathcal{E}_{n-d} + z_{n-d}) \\ &\qquad \wedge \cdots \wedge (\varepsilon_{-(n-d)-1} + \mathcal{E}_1 + z_1) \wedge E_{n-d} \wedge \cdots \wedge E_1\| \\ &= \|(\varepsilon_{-2(n-d)} + \mathcal{E}_{n-d} + z_{n-d}) \\ &\qquad \wedge \cdots \wedge (\varepsilon_{-(n-d)-1} + \mathcal{E}_1 + z_1) \wedge \varepsilon_{n-d+1} \wedge \cdots \wedge \varepsilon_{2(n-d)}\| \cdot |e| \\ &= \|(\varepsilon_{-2(n-d)} + \mathcal{E}_{n-d}) \\ &\qquad \wedge \cdots \wedge (\varepsilon_{-(n-d)} + \mathcal{E}_1) \wedge \varepsilon_{n-d+1} \wedge \cdots \wedge \varepsilon_{2(n-d)}\| \cdot |e| \equiv \alpha(e)|e| \ . \end{aligned} \tag{1.47}$$

From (1.46), we also obtain

$$|\widetilde{m}_{e,z}|^{-1} \geq \|(\varepsilon_{-2(n-d)} + \mathcal{E}_{n-d} + z_{n-d}) \wedge \cdots \wedge (\varepsilon_{-(n-d)-1} + \mathcal{E}_1 + z_1)\| \equiv \alpha(e,z) \ . \tag{1.48}$$

From (1.47) and (1.48)

$$|\widetilde{m}_{e,z}|^{-1} \geq \alpha(e)^{1/2}|e|^{1/2} \cdot \alpha(e,z)^{1/2} \ ;$$

hence, for $t \geq 0$,

$$(|e|\ |\widetilde{m}_{e,z}|)^t \leq (|e|\alpha(e)^{-1}\alpha(e,z)^{-1})^{t/2} \ . \tag{1.49}$$

Let, for a matrix $b$,

$$\lambda_\nu(b) = \begin{cases} \max\{1, \|b\|\} , & \nu < \infty \\ (1+\|b\|^2)^{1/2} , & F_\nu = \mathbb{R} \\ 1+\|b\|^2 , & F_\nu = \mathbb{C} \end{cases} .$$

Then, by the definitions of $\alpha(e)$ and $\alpha(e,z)$,

$$\alpha(e) \geq |e| \quad \text{and} \quad \alpha(e) \geq \lambda_\nu(e)$$
$$\alpha(e,z) \geq |e| \quad \text{and} \quad \alpha(e,z) \geq \lambda_\nu(z) ;$$

hence

$$\alpha(e)\alpha(e,z) \geq |e|\lambda_\nu(e)^{1/2}\lambda_\nu(z)^{1/2} .$$

Using this in (1.49), we get

$$(|e| \; |\widetilde{m}_{e,z}|)^t \leq \lambda_\nu(e)^{-t/4}\lambda_\nu(z)^{-t/4} .$$

Thus the integrability of $A_\nu^{(1)}$ follows from that of

$$\int\limits_{|e|\geq 1} \lambda_\nu(e)^{-t/4} d_a(e) \cdot \int\limits_{\mathcal{U}'_\nu} \lambda_\nu(z)^{-t/4} dz$$

which is valid, for $t \gg 0$, by Proposition 7 in [J-R1]. Similarly,

$$A_\nu^{(2)} = \iint\limits_{|e|>1} (|e^*| \; |m_{e,z}|)^t dedz = \iint\limits_{|e|<1} (|e| \; |\widetilde{m}_{e,z}|)^t \frac{d_a e}{|e|^{n-d}} dz$$
$$\cdot \iint\limits_{|e|<1} |e|^{t-n+d} |\widetilde{m}_{e,z}|^t d_a e dz < \iint\limits_{|e|<1} |\widetilde{m}_{e,z}|^t d_a e dz ,$$

for $t > n-d$. As before

$$|\widetilde{m}_{e,z}|^{-1} \geq \lambda_\nu(e) \quad \text{and} \quad |\widetilde{m}_{e,z}| \geq \lambda_\nu(z) ;$$

hence

$$A_\nu^{(2)} < \iint\limits_{|e|<1} \lambda_\nu(e)^{-t/2}\lambda_\nu(z)^{-t/2} d_a(e) dz \leq \iint\limits_{|e|<1} \lambda_\nu(e)^{-t/4}\lambda_\nu(z)^{-t/4} d_a(e) dz .$$

All in all,

$$A_\nu = A_\nu^{(1)} + A_\nu^{(2)} < \int\limits_{\mathcal{U}'_\nu} \int\limits_{\mathrm{GL}_{n-d}(F_\nu)} \lambda_\nu(e)^{-t/4}\lambda_\nu(z)^{-t/4} d_a(e) dz \equiv r_\nu(t)$$

and the last integral converges absolutely for $t \gg 0$, by Proposition 7 in [J-R1]. By the same proposition, $\prod_\nu r_\nu(t)$ converges for $t \gg 0$. This proves the absolute convergence of $A$ and Lemma 7. □

We go back to the convergence of $I'_{j,d}$ in (1.4.2).

*Case* $j=1,\ 0<d<n$. Here (in (1.42))

$$|g| \leq c(|e^*|\ |m_{e,z}|)^{-1}\ . \tag{1.50}$$

Integral (1.42) splits into two summands $I^{(1)}_{1,d}+I^{(2)}_{1,d}$. In the first summand, we integrate along

$$|g| \leq (|e^*|\ |m_{e,z}|)^{d/n-d}\ . \tag{1.51}$$

We have

$$(|e^*|\ |m_{e,z}|)^{d/n-d} \leq c(|e^*|\ |m_{e,z}|)^{-1}$$

due to Lemma 6 and the fact that $c>1$. Thus (1.50) is redundant in $I^{(1)}_{1,d}$. We have

$$I^{(1)}_{1,d} = \int\limits_{\substack{|g|\leq(|e^*|\ |m_{e,z}|)^{d/n-d}\\ g\in \mathrm{GL}_{2d}(F)\backslash \mathrm{GL}_{2d}(\mathbb{A})\\ e\in \mathrm{GL}_{n-d}(\mathbb{A}),z\in \mathcal{U}'_{\mathbb{A}}}} |g|^{s-n+d-1}(|e^*|\ |m_{e,z}|)^{s+n} dgdedz\ . \tag{1.52}$$

The $dg$-integration in (1.52) is proportional to

$$\int\limits_0^{(|e^*|\ |m_{e,z}|)^{d/n-d}} t^{s-n+d-1}\frac{dt}{t} = \frac{1}{s-n+d-1}(|e^*|\ |m_{e,z}|)^{\frac{d}{n-d}(s-n+d-1)}$$

for $\mathrm{Re}(s)>n-d+1$. $I^{(1)}_{1,d}$ is then proportional to

$$\int (|e^*|\ |m_{e,z}|)^{\frac{d}{n-d}(s-n+d-1)+s+n} dedz \tag{1.53}$$

(that is, it is enough now to get the absolute convergence of (1.53)). By Lemma 7, since $\mathrm{Re}(s)\gg 0$, (1.53) converges absolutely. Now consider $I^{(2)}_{1,d}$.

(1.54)

$$I^{(2)}_{1,d} = \int\limits_{\substack{(|e^*|\ |m_{e,z}|)^{d/n-d}<|g|<c(|e^*|\ |m_{e,z}|)^{-1}\\ g\in \mathrm{GL}_{2d}(F)\backslash \mathrm{GL}_{2d}(\mathbb{A})\\ e\in \mathrm{GL}_{n-d}(\mathbb{A}),\ z\in \mathcal{U}'_{\mathbb{A}}}}$$

$$\cdot\ |g|^{-2(n-d)N+s-n+d-1}(|e^*|\ |m_{e,z}|)^{2dN+s+n} dgdedz\ .$$

The $dg$-integration is proportional to

$$\frac{(c(|e^*|\ |m_{e,z}|)^{-1})^{-(n-d)(2N+1)+s-1}}{-(n-d)(2N+1)+s-1} - \frac{(|e^*|\ |m_{e,z}|)^{\frac{d}{n-d}(-(n-d)(2N+1)+s-1)}}{\frac{d}{n-d}(-(n-d)(2N+1)+s-1)}\ .$$

Thus (1.54) is a difference of two integrals. One is proportional to

$$J_1 = \int (|e^*|\ |m_{e,z}|)^{n(2N+1)+n-d+1} dedz$$

and the other is proportional to

$$J_2 = \int (|e^*|\ |m_{e,z}|)^{\frac{n}{n-d}s+n-d+\frac{d}{n-d}} dedz \ .$$

The integral $J_1$ converges absolutely, since we may choose $N$ as large as we want, and then apply Lemma 7 and similarly, in $J_2$, we assume that $\mathrm{Re}(s) \gg 0$, and get the convergence again.

*Case* $j = 2$ , $0 < d < n$. Here (in (1.42))

$$|g| > c(|e^*|\ |m_{e,z}|)^{-1} \ . \tag{1.55}$$

Again, split (1.42) into two summands $I^{(1)}_{2,d} + I^{(2)}_{2,d}$. In the first summand, integrate along (1.51). Thus, by (1.51) and (1.55),

$$c(|e^*|\ |m_{e,z}|)^{-1} < |g| \le (|e^*|\ |m_{e,z}|)^{d/n-d}$$

and, in particular, $(|e^*|\ |m_{e,z}|)^{n/n-d} > c$. Since $c > 1$, this contradicts Lemma 6. Hence the domain of integration of $I^{(1)}_{2,d}$ is empty, and we need only consider $I^{(2)}_{2,d}$. The integration in $I^{(2)}_{2,d}$ is over $g$ satisfying (1.55) and $|g| > (|e^*|\ |m_{e,z}|)^{d/n-d}$. By Lemma 6, this last condition is redundant (in the presence of (1.55)). Thus,

$$I'_{c,d} = I^{(2)}_{2,d} = \int\limits_{\substack{|g|>c(|e^*|\ |m_{e,z}|)^{-1} \\ g\in \mathrm{GL}_{2d}(F)\backslash \mathrm{GL}_{2d}(\mathbb{A}) \\ e\in \mathrm{GL}_{n-d}(\mathbb{A}),\ z\in \mathcal{U}'_{\mathbb{A}}}} |g|^{-s-(n-d)(2N+1)}(|e^*||\,|m_{e,z}|) dgdedz \ . \tag{1.56}$$

The $dg$-integral in (1.56) is absolutely convergent, once we take $N$ large enough (relative to $s, n, d$), and is proportional to

$$\int\limits_{c(|e^*|\ |m_{e,z}|)^{-1}}^{\infty} t^{-s-(n-d)(2N+1)} \frac{dt}{t} = \frac{(c(|e^*|\ |m_{e,z}|)^{-1})^{-s-(n-d)(2N+1)}}{-s-(n-d)(2N+1)} \ .$$

Thus (1.56) is proportional to

$$\int\limits_{\substack{e\in \mathrm{GL}_{n-d}(\mathbb{A}) \\ z\in \mathcal{U}'_{\mathbb{A}}}} (|e^*|\ |m_{e,z}|)^{2n(N+1)-d+1} dedz \ ,$$

which converges absolutely, provided we choose $N \gg 1$ (Lemma 7).

*Case* $j = 3$, $0 < d < n$. This case follows exactly as Case $j = 2$.

*Case* $j=1,\ d=0$. By (1.38)–(1.40), we have to consider integrals of the form

$$\int\limits_{z\in\mathcal{U}'_{\mathbb{A}}}\ \int\limits_{\substack{e\in\mathrm{GL}_n(\mathbb{A})\\ |e^*|\ |m_{e,z}|<c}}\ \int\limits_{b\in\mathrm{Sp}_{2n}(F)\backslash\mathrm{Sp}_{2n}(\mathbb{A})} \varphi(b\begin{pmatrix} I_n & \\ & e^*\end{pmatrix} m_{e,z})(|e^*|\ |m_{e,z}|)^{s+n}dbdedz$$

for $\varphi$ in the space of $\tau$. Note that, by Lemma 6 (and $c>1$), $|e^*|\ |m_{e,z}|<c$ always holds. Since $\varphi$ is bounded, (and $\mathrm{Sp}_{2n}(F)\backslash\mathrm{Sp}_{2n}(\mathbb{A})$ has finite volume), it is enough to consider the convergence of

$$\int(|e^*|\ |m_{e,z}|)^{s+n}dedz\ .$$

This integral converges absolutely for $\mathrm{Re}(s)$ large enough (Lemma 7).

*Case* $j=2,3;\ d=0$. Similar to case $j=1$, we have to consider

$$\int\limits_{z\in\mathcal{U}'_{\mathbb{A}}}\ \int\limits_{\substack{e\in\mathrm{GL}_n(\mathbb{A})\\ |e^*|\ |m_{e,z}|>c}}\ \int\limits_{b\in\mathrm{Sp}_{2n}(F)\backslash\mathrm{Sp}_{2n}(\mathbb{A})} \varphi(b\begin{pmatrix} I_n & \\ & e^*\end{pmatrix} m_{e,z})(|e^*|\ |m_{e,z}|)^{1-s+n}dbdedz\ .$$

The domain of integration is empty here, by Lemma 6 and the fact that $c>1$.

We proved assertions (1)–(3) of Proposition 5. They provide the justification to the formal proof in Section 4, which shows that $I_{j,d}=0$ for $j=1,2,3$ and $0\le d<n$. This concludes the proof of Theorem 2. □

## 2. A tower of correspondences for $\mathrm{GL}_m$

1. *Definition of the tower.* Let $\tau$ be an irreducible, automorphic, cuspidal, representation of $\mathrm{GL}_{2n}(\mathbb{A})$. Assume that $\tau$ is self-dual. As recalled in the introduction, exactly one of the partial $L$-functions $L^S(\tau,\Lambda^2,s)$ or $L^S(\tau,\mathrm{Sym}^2,s)$ has a pole (simple) at $s=1$. Consider the Eisenstein series on $\mathrm{SO}_{4n}(\mathbb{A})$, in the first case, and on $\mathrm{SO}_{4n+1}(\mathbb{A})$, in the second case, induced from the Siegel-type parabolic subgroup, the representation $\tau$ and the same complex parameter $s$, as in Chapter 1 Denote it by $E(g,f^{\phi}_{\tau,s})$, keeping the notation of Chapter 1. In subsection 1.1, we have seen that it has a simple pole at $s=1$. Let

$$E_1(g,\phi)=\operatorname*{Res}_{s=1}E(g,f^{\phi}_{\tau,s})\ .$$

By taking a sequence of certain Fourier coefficients of $E_1$ along certain unipotent subgroups of $\mathrm{SO}_{4n}$, in the first case, and of $\mathrm{SO}_{4n+1}$, in the second case, and restricting to subgroups of the form $\mathrm{SO}_{2k+1}(\mathbb{A})$ (resp. $\mathrm{SO}_{2k}(\mathbb{A})$), we obtain a "tower" of spaces of automorphic representations $\sigma_k(\tau)$ of $\mathrm{SO}_{2k+1}(\mathbb{A})$ (resp. of $\mathrm{SO}_{2k}(\mathbb{A})$) which are generic

$$\begin{array}{cccccc} & & & & \mathrm{SO}_{4n-1}(\mathbb{A}) \\ & & & \nearrow & \vdots \\ \tau & \text{on} & \mathrm{GL}_{2n}(\mathbb{A}) & \longrightarrow & \mathrm{SO}_{2n+1}(\mathbb{A}) \\ & & & \searrow & \vdots \\ & & & \searrow & \mathrm{SO}_3(\mathbb{A}) \\ & & & & \mathrm{SO}_1(\mathbb{A}) \end{array}, \tag{2.1}$$

$$\begin{array}{cccccc} & & & & \mathrm{SO}_{4n}(\mathbb{A}) \\ & & & \nearrow & \vdots \\ \tau & \text{on} & \mathrm{GL}_{2n}(\mathbb{A}) & \longrightarrow & \mathrm{SO}_{2n}(\mathbb{A}) \\ & & & \searrow & \vdots \\ & & & \searrow & \mathrm{SO}_2(\mathbb{A}) \\ & & & & \{I\} \end{array}. \tag{2.2}$$

We will prove that the last step in the tower is always nontrivial and that in the first step $k$, where $\sigma_k(\tau)$ is nonzero, $\sigma_k(\tau)$ is also cuspidal. In this sense, the towers resemble Rallis original symplectic-orthogonal towers. In case (2.1), we can be more precise and distinguish two cases, according to whether $L^S(\tau, st, \frac{1}{2}) \neq 0$ or not. If $L^S(\tau, st, \frac{1}{2}) \neq 0$, consider the Eisenstein series $E(g, f^{\phi}_{\tau,s})$ on $\mathrm{Sp}_{4n}(\mathbb{A})$ (as in Chapter 1). Again, by 1.1, $E_1(g,\phi) = \operatorname*{Res}_{s=1} E(g, f^{\phi}_{\tau,s})$ is nontrivial, and by taking a sequence of certain Fourier-Jacobi coefficients of $E_1(g,\phi)$, we construct another tower of liftings to automorphic, generic representations of $\widetilde{\mathrm{Sp}}_{2k}(\mathbb{A})$

$$\begin{array}{cccccc} & & & & \widetilde{\mathrm{Sp}}_{4n-2}(\mathbb{A}) \\ & & & \nearrow & \vdots \\ \tau & \text{on} & \mathrm{GL}_{2n}(\mathbb{A}) & \longrightarrow & \widetilde{\mathrm{Sp}}_{2n}(\mathbb{A}) \\ & & & \searrow & \vdots \\ & & & & \widetilde{\mathrm{Sp}}_2(\mathbb{A}) = \widetilde{\mathrm{SL}}_2(\mathbb{A}) \\ & & & \searrow & \{I\} \end{array}. \tag{2.3}$$

*Towers* (2.1) *and* (2.2). Let $G$ be one of the groups $\mathrm{SO}_{4n}$ or $\mathrm{SO}_{4n+1}$, acting from the left on the column space $W = F^{4n}$ (resp. $W = F^{4n+1}$) equipped with the quadratic form defined by $\begin{pmatrix} & & 1 \\ & \iddots & \\ 1 & & \end{pmatrix}$. Let $\{\varepsilon_1, \ldots, \varepsilon_{2n}, \varepsilon_{-2n}, \ldots, \varepsilon_{-1}\}$ (resp. $\{\varepsilon_1, \ldots, \varepsilon_{2n}, \varepsilon_0, \varepsilon_{-2n}, \ldots, \varepsilon_{-1}\}$) be the standard basis of $W$. Let $V_\ell^+ = \mathrm{Span}\{\varepsilon_1, \ldots, \varepsilon_\ell\}$ $(1 \le \ell \le 2n)$. This is an $\ell$-dimensional isotropic subspace of $W$. Its dual in $W$ is $V_\ell^- = \mathrm{Span}\{\varepsilon_{-1}, \ldots, \varepsilon_{-\ell}\}$. Consider the parabolic subgroup $P_k$ of $G$, which preserves $E_k = V_{2n-k-1}^+$ (resp. $E_k = V_{2n-k}^+$). Its Levi decomposition is $P_k = M_k \ltimes U_k$. $U_k$ is abelian, in case $\dim E_k = 1$, or in case $\dim E_k = 2n$, when $G = \mathrm{SO}_{4n}$, and otherwise, it is a two-step nilpotent group. The Siegel parabolic subgroup $P$ corresponds to $k = -1$ (resp. $k = 0$) Let $E_k'$ be the dual of $E_k$ inside $W$, and let $W_k$ be the orthocomplement of $E_k + E_k'$ in $W$. Then

$$M_k \cong \mathrm{GL}(E_k) \times \mathrm{SO}(W_k) \ ,$$

and $U_k^{ab}$ is isomorphic to the space of matrices of size $\dim E_k \times \dim W_k$. Let $\alpha_k = \dim E_k$, $\beta_k = \dim W_k$. In matrix form, an element of $M_k$ has the form $m = \begin{pmatrix} a & & \\ & b & \\ & & a^* \end{pmatrix}$, where $a \in \mathrm{GL}(E_k)$ and $b \in \mathrm{SO}(W_k)$, and an element of $U_k^{ab}$ has the form

$$u = \begin{pmatrix} I_{\alpha_k} & Y & * \\ & I_{\beta_k} & Y' \\ & & I_{\alpha_k} \end{pmatrix} \ . \tag{2.4}$$

Conjugation by $m$ gives

$$mum^{-1} = \begin{pmatrix} I_{\alpha_k} & aYb^{-1} & * \\ & I_{\beta_k} & bY'(a^*)^{-1} \\ & & I_{\alpha_k} \end{pmatrix} \ . \tag{2.5}$$

If we realize $U_k^{ab}$ as $E_k \otimes W_k^*$, where $W_k^*$ is the dual of $W_k$, written in the form of rows in $\beta_k$ coordinates, then conjugation (2.5) by $m$ of $u$, which corresponds to $e \otimes v$, gives $ae \otimes v \cdot b^{-1}$. Define

$$v_0 = \begin{cases} {}^t\varepsilon_0, & \dim W = 4n+1 \\ {}^t\varepsilon_{2n} + {}^t\varepsilon_{-2n}, & \dim W = 4n \ . \end{cases}$$

Denote (again) the symmetric bilinear form on $W_k^*$ by $(,)$, where $({}^tv_1, {}^tv_2) = (v_1, v_2)$, for $v_1$ and $v_2$ in $W_k$. Denote

$$y_k = \varepsilon_{-\alpha_k} \otimes v_0 \ . \tag{2.6}$$

Fix a nontrivial character $\psi$ of $F \backslash \mathbb{A}$ and define the following character $\psi_k$ of $U_k(F) \backslash U_k(\mathbb{A})$. Let $u$ in $U_k^{ab}(\mathbb{A})$ correspond to $y \in E_k(\mathbb{A}) \otimes W_k^*(\mathbb{A})$. Then

$$\psi_k(u) = \psi\big((y, y_k)\big) \ , \tag{2.7}$$

(where ( , ) is extended naturally to $(E_k + E_k')(\mathbb{A}) \otimes W_k^*(\mathbb{A})$). In matrix form ($u$ as in (2.4))

$$\psi_k(u) = \begin{cases} \psi(Y_{2n-k,k+1}), & \dim W = 4n+1 \\ \psi(Y_{2n-k-1,k+1} + Y_{2n-k-1,k+2}), & \dim W = 4n \end{cases} \ . \tag{2.8}$$

The stabilizer of $\psi_k$ in $\mathrm{SO}(W_k)$ is $H_k = \mathrm{SO}(W'_k)$, where $W'_k$ is the orthocomplement of ${}^t v_0$ in $W_k$. Note that

$$(2.9) \qquad H_k \cong \begin{cases} \mathrm{SO}_{2k}, & \dim W = 4n+1 \\ \mathrm{SO}_{2k+1}, & \dim W = 4n \end{cases} .$$

The stabilizer of $\psi_k$ in $\mathrm{GL}(E_k)$ is the mirabolic subgroup (of type $(\alpha_k - 1, 1)$); in particular, it contains $Z_{\alpha_k}$ – the standard maximal unipotent subgroup of $\mathrm{GL}(E_k)$. Extend $\psi_k$ to $Z_{\alpha_k}(\mathbb{A})$ by the standard nondegenerate character of $Z_{\alpha_k}(\mathbb{A})$, defined by $\psi$, i.e. it takes $z$ in $Z_{\alpha_k}(\mathbb{A})$ to $\psi(\sum_{\ell=1}^{\alpha_k-1} z_{\ell,\ell+1})$. This defines a character $\psi_k$ of $N_k(\mathbb{A})$, where $N_k = Z_{\alpha_k} \cdot U_k$.

Consider the Eisenstein series $E(g, f^{\phi}_{\tau,s})$ of Chapter 1 (mentioned in Remark 2, following Proposition 1). Recall that $\tau$ is an irreducible, automorphic cuspidal representation of $\mathrm{GL}_{2n}(\mathbb{A})$, which is self-dual. $\phi$ is a smooth $K$-finite function in $\mathrm{Ind}_{P(\mathbb{A})}^{G(\mathbb{A})} \tau$. $P = M \ltimes U$ is the Siegel-type parabolic subgroup of $G = \mathrm{SO}(W)$. Finally

$$(2.10) \qquad \begin{aligned} \varphi^{\phi}_{\tau,s}(g;m) &= H(g)^{s-1/2}\phi(g;m), \qquad g \in G_{\mathbb{A}}, \qquad m \in \mathrm{GL}_{2n}(\mathbb{A}) \\ f^{\phi}_{\tau,s}(g) &= \varphi^{\phi}_{\tau,s}(g;1) \ , \end{aligned}$$

and if $g$ is written in the Iwasawa decomposition $\widehat{a}uk$, for $a \in \mathrm{GL}_{2n}(\mathbb{A})$, $u \in U(\mathbb{A})$, $k \in$K, then $H(g) = |\det a|$.

Assume that $E(g, f^{\phi}_{\tau,s})$ has a pole at $s = 1$. Let $\sigma_k(\tau)$ be the representation, by right translations, of $H_k(\mathbb{A})$ on the space spanned by the $\psi_k$-Fourier coefficients of $E_1(g,\phi)$, along $N_k$, restricted to $H_k(\mathbb{A})$. That is, the space of $\sigma_k(\tau)$ is the space of automorphic functions on $H_k(\mathbb{A})$, of the form (2.11)

$$(2.11) \qquad p_k(h) = p_k(h,\phi) = \int\limits_{N_k(F)\backslash N_k(\mathbb{A})} \overline{E_1(vh,\phi)}\psi_k(v)dv \ .$$

*The tower* (2.3). Here $G = \mathrm{Sp}_{4n}$. Let $W = F^{4n}$ (column space), with the symplectic form ( , ) defined by $J_{4n} = \begin{pmatrix} &&&&&1 \\ &&&&\iddots& \\ &&&1&& \\ &&-1&&& \\ &\iddots&&&& \\ -1&&&&& \end{pmatrix}$. Let $\{\varepsilon_1, \dots, \varepsilon_{2n}, \varepsilon_{-2n}, \dots, \varepsilon_{-1}\}$ be the standard basis of $W$. Define $V^{\pm}_{\ell}$, $1 \le \ell \le 2n$, as in the previous case, and $E_k = V^{+}_{2n-k}$. Let $P_k$ be the parabolic subgroup of $G$, which preserves $E_k$. Its Levi decomposition is $P_k = M_k \ltimes U_k$. $U_k$ is abelian in case $\dim E_k = 2n$ (i.e. $k = 0$), and otherwise it is a two step nilpotent group. Let $E'_k$ and $W_k$ be as in the previous case. Then $M_k \cong \mathrm{GL}(E_k) \times \mathrm{Sp}(W_k) = \mathrm{GL}_{2n-k} \times \mathrm{Sp}_{2k}$. The Siegel parabolic subgroup $P$ corresponds to $k = 0$. Let $N_k = Z_{2n-k} \cdot U_k$ ($Z_{2n-k}$ is the standard maximal unipotent subgroup of $\mathrm{GL}(E_k) = \mathrm{GL}_{2n-k}$). Define a character $\psi_k$ on $N_k(\mathbb{A})$ to be trivial on $U_k(\mathbb{A})$ and the standard nondegenerate character (corresponding to $\psi$) on $Z_{2n-k}(\mathbb{A})$. Note that $\psi_k$ is fixed by $\mathrm{Sp}(W_k)_{\mathbb{A}}$. The group $N_k/N_{k+1}$ is isomorphic to the Heisenberg group $\mathcal{H}_k$ in $2k+1$ variables. In matrix form, $v$ in $N_k$ looks like this

$$v = \begin{pmatrix} z & * & * & * & * \\ & 1 & e & y & * \\ & & I_{2k} & e' & * \\ & & & 1 & * \\ & & & & z^* \end{pmatrix}, \quad z \in Z_{2n-(k+1)} \ .$$

Then the isomorphism takes $v$ modulo $N_{k+1}$ to $(\tilde{e}; y)$ where $\tilde{e} = e \cdot \begin{pmatrix} I_k & \\ & 2I_k \end{pmatrix}$. Denote the composition of this isomorphism with the quotient map $N_k \to N_k/N_{k+1}$ by $j_k$. Consider the Weil representation of $\widetilde{\mathrm{Sp}}_{2k}(\mathbb{A})$ – the global metaplectic cover. It comes with a choice of a nontrivial additive character of $F\backslash\mathbb{A}$. We choose $\psi$. Denote the Weil representation by $\omega_{\psi,k}$, and extend it to $\mathcal{H}_k(\mathbb{A}) \cdot \widetilde{\mathrm{Sp}}_{2k}(\mathbb{A})$ by the Stone von-Neumann representation (of central characters $\psi$) on $\mathcal{H}_k(\mathbb{A})$. We still denote this representation by $\omega_{\psi,k}$. It acts on the space $S(\mathbb{A}^k)$ of Schwartz-Bruhat functions on $\mathbb{A}^k$. Consider its automorphic realization on the space of theta series

$$(2.12) \qquad \theta^{\xi}_{\psi,k}(h) = \sum_{x \in F^k} \omega_{\psi,k}(h)\xi(x); \qquad h \in \mathcal{H}_{\mathbb{A}} \cdot \widetilde{\mathrm{Sp}}_{2k}(\mathbb{A}), \qquad \xi \in S(\mathbb{A}^k) \ .$$

Let $\tau$ be an irreducible, automorphic, cuspidal representation of $\mathrm{GL}_{2n}(\mathbb{A})$. Assume that $L^S(\tau, \Lambda^2, s)$ has a pole at $s = 1$ and that $L(\tau, st, \frac{1}{2}) \neq 0$. Then the Eisenstein series $E(g, f^{\phi}_{\tau,s})$ has a (simple) pole at $s = 1$ (subsection 1.1). Let $\sigma_k(\tau)$ be the representation, by right translations, of $\widetilde{\mathrm{Sp}}_{2k}(\mathbb{A})$, on the space spanned by coefficients of Fourier-Jacobi type on $E_1(g, \phi)$, namely the functions $p_k(h)$, where

$$(2.13) \qquad p_k(h) = p_k(h, \phi, \xi) = \int\limits_{N_k(F)\backslash N_k(\mathbb{A})} E_1(vh, \phi)\theta^{\xi}_{\psi,k}(j_k(v)h)\psi_k^{-1}(v)dv \ .$$

Here $h \in \widetilde{\mathrm{Sp}}_{2k}(\mathbb{A})$, $\phi$- as in subsection 1.1 and $\xi \in S(\mathbb{A}^k)$. The appearance of $h$ in $E_1(vh, \phi)$ factors through $\mathrm{Sp}_{2k}(\mathbb{A})$. Denote, for this tower $H_k = \widetilde{\mathrm{Sp}}_{2k}$.

*The case of* $\mathrm{GL}_m$, *where* $m$ *is odd.* Assume that $m = 2n+1$. Let $\tau$ be an irreducible, automorphic, cuspidal representation of $\mathrm{GL}_{2n+1}(\mathbb{A})$. Assume that $\tau$ is self-dual. In this case, it is known that $L^S(\tau, \Lambda^2, s)$ is holomorphic, and hence $L^S(\tau, \mathrm{Sym}^2, s)$ has a pole at $s = 1$. Consider the Eisenstein series $E(g, f^{\phi}_{\tau,s})$ on $G_{\mathbb{A}} = \widetilde{\mathrm{Sp}}_{4n+2}(\mathbb{A})$. It has a pole at $s = 1$. As in the last case, consider the groups $U_k$ and $N_k$ (only now $\dim E_k = 2n - k + 1$). $\psi_k$ is the character of $N_k(\mathbb{A}) = Z_{2n-k+1}(\mathbb{A})U_k(\mathbb{A})$, which is trivial on $U_k(\mathbb{A})$ and is the nondegenerate character on $Z_{2n-k+1}(\mathbb{A})$. As in the tower (2.3) and (2.13), let $\sigma_k(\tau)$ be the representation, by right translations, of $\mathrm{Sp}_{2k}(\mathbb{A})$ on the space spanned by coefficients of Fourier-Jacobi type on $\overline{E_1(g, \phi)}$ namely the functions $\overline{p_k(h)}$, where

$$p_k(h) = p_k(h, \phi, \xi) = \int\limits_{N_k(F)\backslash N_k(\mathbb{A})} E_1(vh, \phi)\theta^{\xi}_{\psi,k}(j_k(v)h)\psi_k^{-1}(v)dv \ . \tag{2.14}$$

Here $h \in \mathrm{Sp}_{2k}(\mathbb{A})$, $\phi$ – as in subsection 1.1 and $\xi \in S(\mathbb{A}^k)$. Note that in (2.14), we may take any preimage of $h$ in $\widetilde{\mathrm{Sp}}_{2k}(\mathbb{A})$.

This defines a tower for $\tau$ on $\mathrm{GL}_{2n+1}(\mathbb{A})$

$$\begin{array}{ccccc} & & & & \mathrm{Sp}_{4n}(\mathbb{A}) \\ & & & \nearrow & \vdots \\ \tau & \text{on} & \mathrm{GL}_{2n+1}(\mathbb{A}) & \longrightarrow & \mathrm{Sp}_{2n}(\mathbb{A}) \\ & & & \searrow & \vdots \\ & & & & \mathrm{Sp}_2(\mathbb{A}) \\ & & & \searrow & \{I\} \end{array} \ . \tag{2.15}$$

Denote, for this tower, $H_k = \mathrm{Sp}_{2k}$.

The definition of $\sigma_k(\tau)$ through formulae (2.11), (2.13) and (2.14) is not accidental. It is prescribed by Rankin-Selberg integrals which represent the standard $L$-function for $H_k \times \mathrm{GL}_m$ (see [G-PS], [G],[So], [GRS2]). More precisely, let $\sigma$ be an irreducible, automorphic, cuspidal representation of $H_k(\mathbb{A})$ ($H_k$ is considered in each of the cases $G = \mathrm{SO}_{4n}, \mathrm{SO}_{4n+1}, \mathrm{Sp}_{4n}, \widetilde{\mathrm{Sp}}_{4n+2}$). Let $\varphi$ be a cusp form of $\sigma$, then the integral

$$\int\limits_{H_k(F)\backslash H_k(\mathbb{A})} \varphi(h)\overline{p_k(h)}dh \tag{2.16}$$

is the residue at $s = 1$ of the Rankin-Selberg integrals mentioned above (it is obtained from (2.16), by replacing $E_1(h, \phi)$ with $E(h, f^{\phi}_{\tau,s})$ in the definition of $p_k(h)$.) We remark that similar towers of unitary groups can be constructed with the aid of similar global integrals (the theory should resemble the previous cases; see [T], [Wt].)

2. *The tower of correspondences.* Let $G$ be one of the groups $\mathrm{SO}_{4n}$, $\mathrm{SO}_{4n+1}, \mathrm{Sp}_{4n}$, $\widetilde{\mathrm{Sp}}_{4n+2}$, and let $\tau$ be an irreducible, automorphic, cuspidal and self-dual representation of $\mathrm{GL}_m(\mathbb{A})$, where $m = 2n$ or $m = 2n + 1$, such that the corresponding Eisenstein series $E(g, f^{\phi}_{\tau,s})$ has a pole at $s = 1$. Consider the representation $\sigma_k(\tau)$ of $H_k(\mathbb{A})$, acting by right translations in the space spanned by the functions $p_k(h)$ ((2.11), (2.13), (214)).

Our main theorem in this chapter is:

THEOREM 8. *Assume that $\sigma_k(\tau) = 0$, for $1 \leq k \leq \ell - 1$. Then $\sigma_\ell(\tau)$ is either zero or cuspidal on $H_\ell(\mathbb{A})$.*

This is the well-known property of the symplectic-orthogonal towers of theta liftings [R]. (See also [GRS3], for another tower for $G_2$, within exceptional groups.) The towers we present here "should start" at level $k = n$. We cannot prove this at this stage(*), but at the third chapter, we will prove that, in case $G = \mathrm{Sp}_{4n}$, $\sigma_k(\tau) = 0$, for $1 \leq k < n$. Once we prove that $\sigma_n(\tau)$ is nontrivial, then $\sigma_n(\tau)$ is a cuspidal representation of $H_n(\mathbb{A})$, associated to $\tau$ (such that $E(g, f^{\phi}_{\tau,s})$ has a pole at $s = 1$). The correspondence $\tau \mapsto \sigma_n(\tau)$ should be the explicit backward map to the functorial lift from $\mathrm{GL}_m$ to $H_n$ ($m = 2n, 2n+1$).

We first prove (recalling that $P = M \ltimes U$ denotes the Siegel parabolic subgroup of $G$),

PROPOSITION 9.

1. *The constant term of $E_1(g, \phi)$ along any unipotent radical of a standard parabolic subgroup, different from $P$, is zero.*
2. *The Whittaker coefficient of $E_1(g, \phi)$, along the maximal unipotent subgroup of $G$, is zero.*
3. *Let $Z_m$ be the standard maximal unipotent subgroup of $\mathrm{GL}_m$, embedded naturally in $M$ — the Levi part of $P$. Consider the standard nondegenerate character of $Z_m(\mathbb{A})$, defined by $\psi$, and continue to denote it by $\psi$. Then the following coefficient is nontrivial,*

$$\int\limits_{Z_m(F)\backslash Z_m(\mathbb{A})} \int\limits_{U_F\backslash U_{\mathbb{A}}} E_1(zug)\psi(z)dzdu \ .$$

*Proof.* Part 1 follows, since $E(g, f^{\phi}_{\tau,s})$ is concentrated on $P$. Part 2 follows since the Whittaker coefficient of $E(g, f^{\phi}_{\tau,s})$ is holomorphic at $s = 1$. (The proof uses the formula on p. 352 of [Sh], for the Whittaker coefficient of $E(g, f^{\phi}_{\tau,s})$. The product of $L$-functions, which appears there, is that given by the denominator in each case of (1.5). As explained (in parenthesis) in the proof of Proposition 1, each such denominator is nonzero (and holomorphic) at $s = 1$. To complete the argument, use Proposition 3.1 in [Sh]. The proof in case $E(g, f^{\phi}_{\tau,s})$ is on $\widetilde{\mathrm{Sp}}_{4n}$ is exactly the same.) As for part 3, consider the constant term of $E_1(g,\phi)$ along $U$. Then by part 1, it defines a cusp form $\xi$ on $\mathrm{GL}_m(\mathbb{A}) \cong M_{\mathbb{A}}$ and so $\xi$ has a nontrivial Whittaker coefficient (with respect to $Z_m$ and $\psi$). This proves Proposition 9. □

(*) *Added in proof*. As we mentioned in the introduction, we now can prove that $\sigma_n(\tau) \neq 0$. The case $G = \mathrm{Sp}_{4n}$ appears in [G-R-S5]. The remaining cases will appear in a work under preparation. The fact that the unramified parameters of $\tau$ and $\sigma_n(\tau)$ correspond, for $G = \mathrm{Sp}_{4n}$ appears in [G-R-S6].

*Proof of Theorem* 8. The proof uses Fourier expansions and Proposition 9 (parts 1 and 2). We present the proof in cases $G = \mathrm{SO}_{4n+1}$ and $G = \mathrm{Sp}_{4n}$. The proof in the two remaining cases is entirely similar.

*Case* $G = \mathrm{SO}_{4n+1}$. Let $R_p$ denote the unipotent radical of the standard, maximal parabolic subgroup of $H_\ell = \mathrm{SO}_{2\ell}$, whose Levi part is isomorphic to $\mathrm{GL}_p \times \mathrm{SO}_{2(\ell-p)}$ $(1 \le p < \ell)$. The elements of $R_p$, naturally identified with matrices in $\mathrm{SO}_{2\ell}$, are of the form

$$\begin{pmatrix} I_p & X & Y \\ & I_{2(\ell-p)} & X' \\ & & I_p \end{pmatrix} . \tag{2.17}$$

The corresponding element in $\mathrm{SO}_{4n+1}$ is

$$r = \begin{pmatrix} I_{2n-\ell} & & & & & & \\ & I_p & X_1 & 0 & X_2 & Y & \\ & & I_{\ell-p} & 0 & 0 & X_2' & \\ & & & 1 & 0 & 0 & \\ & & & & I_{\ell-p} & X_1' & \\ & & & & & I_p & \\ & & & & & & I_{2n-\ell} \end{pmatrix} , \tag{2.18}$$

where $X = (X_1, X_2)$. We want to prove that

$$\int\limits_{R_p(F)\backslash R_p(\mathbb{A})} p_\ell(r,\phi)dr \equiv 0, \qquad \text{for all} \quad 1 \le p < \ell$$

which, by (2.11), means that

$$\int\limits_{R_p(F)\backslash R_p(\mathbb{A})} \int\limits_{N_\ell(F)\backslash N_\ell(\mathbb{A})} E_1(vr,\phi)\psi_\ell(v)dvdr \equiv 0 \ . \tag{2.19}$$

We want to relate the integral (2.19) to an integral over the group $N_{\ell-p}$. For this, it is convenient to conjugate the variables in (2.19) by the Weyl element

$$\beta = \begin{pmatrix} & I_p & & & \\ I_{2n-\ell} & & & & \\ & & I_{2(\ell-p)+1} & & \\ & & & & I_{2n-\ell} \\ & & & I_p & \end{pmatrix} .$$

Conjugating $vr$, in (2.19), by $\beta$, we get a matrix of the form

$$\beta vr\beta^{-1} = \begin{pmatrix} I_p & 0 & X_1 & t & X_2 & D & Y \\ L & Z & A & b & C & E & D' \\ & & I_{\ell-p} & 0 & 0 & C' & X_2' \\ & & & 1 & 0 & b' & t' \\ & & & & I_{\ell-p} & A' & X_1' \\ & & & & & Z^* & 0 \\ & & & & & L' & I_p \end{pmatrix} , \tag{2.20}$$

with $t = 0$ and $Z \in Z_{2n-\ell}$. Note that

$$\psi_\ell(v) = \psi\left(\sum_{j=1}^{2n-\ell-1} Z_{j,j+1}\right)\psi(b_{2n-\ell}) \ .$$

Denote $S = \beta N_\ell R_p \beta^{-1}$. The identity (2.19) (after replacing $\beta \cdot \phi$ by $\phi$) is equivalent to

$$\int\limits_{S_F\backslash S_{\mathbb{A}}} E_1(s,\phi)\psi_\ell(\beta^{-1}s\beta)ds \equiv 0 \ . \tag{2.21}$$

Let $\widetilde{S}$ be the subgroup of matrices (2.20) with $L = 0$, $t = 0$ and $Z = I_{2n-\ell}$. Denote by $\mathcal{L}$ the subgroup of matrices (2.20), that have all coordinates above the diagonal equal to zero. As usual, we identify $Z_{2n-\ell}$ as a subgroup of $S$. Then the left-hand side of (2.21) equals

$$\int\limits_{Z_{2n-\ell}(F)\backslash Z_{2n-\ell}(\mathbb{A})} \int\limits_{\mathcal{L}_F\backslash\mathcal{L}_{\mathbb{A}}} \int\limits_{\widetilde{S}_F\backslash\widetilde{S}_{\mathbb{A}}} E_1(sxz,\phi)\psi_\ell(\beta^{-1}sz\beta)dsdxdz \ . \tag{2.22}$$

Denote by $u_t$ a matrix of the form (2.20), with all other variables, except for $t$ (and any appropriate $Y$) being zero. Consider the function of $F^p\backslash\mathbb{A}^p$, $t \mapsto \int\limits_{\widetilde{S}_F\backslash\widetilde{S}_{\mathbb{A}}} E_1(u_t sh,\phi)\psi_\ell(\beta^{-1}s\beta)ds$ ($h$ is fixed). Take its Fourier expansion and substitute $t = 0$. Then it is easy to see that it equals

$$\sum_{\substack{\eta=(\eta_1,\dots,\eta_p)\\ \eta_j\in F}} \int\limits_{U_{\ell-p}(F)\backslash U_{\ell-p}(\mathbb{A})} E_1(u\eta^{(0)}h,\phi)\psi_{\ell-p}(u)du \ , \tag{2.23}$$

where

$$\eta^{(0)} = \begin{pmatrix} I_p & & & & & & \\ 0 & I_{2n-\ell-1} & & & & & \\ \eta & 0 & 1 & & & & \\ & & & I_{2(\ell-p)+1} & & & \\ & & & & 1 & & \\ & & & & & I_{2n-\ell-1} & \\ & & & & \eta' & 0 & I_p \end{pmatrix} .$$

Recall that $U_k$ is the unipotent radical of the standard parabolic subgroup $P_k$. Using (2.23), we see that (2.22) equals

$$\int\limits_{Z_{2n-\ell}(F)\backslash Z_{2n-\ell}(\mathbb{A})} \int\limits_{\mathcal{L}'_F\backslash\mathcal{L}_{\mathbb{A}}} \int\limits_{U_{\ell-p}(F)\backslash U_{\ell-p}(\mathbb{A})} E_1(uxz,\phi)\psi_{\ell-p}(uz)dudxdz \ , \tag{2.24}$$

where $\mathcal{L}'$ is the subgroup of matrices (2.20) in $\mathcal{L}$, such that the last line of $L$ is zero. Denote by $\mathcal{L}^{(0)}$ the subgroup of matrices (2.20) in $\mathcal{L}$ such that $L$ has all lines zero except for the last one. Write

$$\mathcal{L}' = \mathcal{L}^{(1)}\widetilde{\mathcal{L}}$$
$$Z_{2n-\ell} = Z^{(1)}\widetilde{Z}$$

where $\mathcal{L}^{(1)}$ consists of all matrices (2.20) in $\mathcal{L}'$ such that $L$ has all lines zero except for line $2n-\ell-1$, while $\widetilde{\mathcal{L}}$ corresponds to $L$ with line $2n-\ell-1$ (and of course line $2n-\ell$) being zero. Also

$$Z^{(1)} = \left\{ \begin{pmatrix} I_{2n-\ell-1} & * \\ & 1 \end{pmatrix} \right\} \quad \text{and} \quad \widetilde{Z} = \left\{ \begin{pmatrix} \zeta & 0 \\ & 1 \end{pmatrix} \in Z_{2n-\ell} \right\} .$$

Then (2.24) equals

$$\int\limits_{\mathcal{L}^{(0)}_{\mathbb{A}}} \int\limits_{\widetilde{Z}_F\backslash\widetilde{Z}_{\mathbb{A}}} \int\limits_{\widetilde{\mathcal{L}}_F\backslash\widetilde{\mathcal{L}}_{\mathbb{A}}} \int\limits_{Z^{(1)}_F\backslash Z^{(1)}_{\mathbb{A}}} \int\limits_{\mathcal{L}^{(1)}_F\backslash\mathcal{L}^{(1)}_{\mathbb{A}}} \int\limits_{U_{\ell-p}(F)\backslash U_{\ell-p}(\mathbb{A})} \cdot\ E_1(ux^{(1)}z^{(1)}\widetilde{x}\widetilde{z}x^{(0)},\phi)\psi_{\ell-p}(uz^{(1)}\widetilde{z})d(\cdots) \ . \tag{2.25}$$

Note that $x^{(1)} \in \mathcal{L}^{(1)}_{\mathbb{A}}$ and $z^{(1)} \in Z^{(1)}_{\mathbb{A}}$ commute. Consider the function on $F^p\backslash\mathbb{A}^p$, $a \mapsto \int\limits_{U_{\ell-p}(F)\backslash U_{\ell-p}(\mathbb{A})} E_1(u\widetilde{a}h)\psi_{\ell-p}(u)du$ ($h$ fixed), where

$$\widetilde{a} = \begin{pmatrix} I_p & 0 & a &&&& \\ & I_{2n-\ell-1} & 0 &&&& \\ && 1 &&&& \\ &&& I_{2(\ell-p)+1} &&& \\ &&&& 1 & 0 & a' \\ &&&&& I_{2n-\ell-1} & 0 \\ &&&&&& I_p \end{pmatrix} .$$

The Fourier expansion of this function (at $a=0$) is

$$\sum_{\eta\in F^p} \int\limits_{F^p\backslash\mathbb{A}^p} \int\limits_{U_{\ell-p}(F)\backslash U_{\ell-p}(\mathbb{A})} E_1(u\widetilde{a}h)\psi_{\ell-p}(u)\psi(\eta\cdot a)duda \ . \tag{2.26}$$

Let

$$\eta^{(1)} = \begin{pmatrix} I_p &&&&&& \\ 0 & I_{2n-\ell-2} &&&&& \\ \eta & 0 & 1 &&&& \\ &&& I_{2(\ell-p)+3} &&& \\ &&&& 1 & I_{2n-\ell-2} & \\ &&&& \eta' & 0 & I_p \end{pmatrix} .$$

Then, fixing $\widetilde{x}, \widetilde{z}, x^{(0)}$ in (2.25), the inner-triple integration equals, using (2.26),

$$\int\limits_{\mathcal{L}_F^{(1)}\backslash\mathcal{L}_{\mathbb{A}}^{(1)}} \sum_{\eta\in F^p} \int\limits_{Z_F^{(1)}\backslash Z_{\mathbb{A}}^{(1)}} \int\limits_{F^p\backslash\mathbb{A}^p} \int\limits_{U_{\ell-p}(F)\backslash U_{\ell-p}(\mathbb{A})} \cdot E_1\big(u\cdot(\eta^{(1)}\widetilde{a}\eta^{(1)^{-1}})z^{(1)}$$

$$\cdot\,\eta^{(1)}x^{(1)}\cdot\widetilde{x}\widetilde{z}x^{(0)},\phi\big)\cdot\psi(\eta\cdot a)\psi_{\ell-p}(uz^{(1)})dudadz^{(1)}dx^{(1)} .$$

Note that $\eta^{(1)}\widetilde{a}\eta^{(1)^{-1}}z^{(1)} = \widetilde{a}z'^{(1)}$ and $\psi(\eta\cdot a)\psi_{\ell-p}(z^{(1)}) = \psi_{\ell-p}(z'^{(1)})$. We get that (2.25) equals

$$\text{(2.27)}\qquad \int\limits_{\mathcal{L}_{\mathbb{A}}^{(0)}\mathcal{L}_{\mathbb{A}}^{(1)}} \int\limits_{\widetilde{Z}_F\backslash\widetilde{Z}_{\mathbb{A}}} \int\limits_{\widetilde{\mathcal{L}}_F\backslash\widetilde{\mathcal{L}}_{\mathbb{A}}} \int\limits_{Z_F^{(1)}\backslash Z_{\mathbb{A}}^{(1)}} \int\limits_{F^p\backslash\mathbb{A}^p} \int\limits_{U_{\ell-p}(F)\backslash U_{\ell-p}(\mathbb{A})}$$

$$\cdot\ E_1(u\widetilde{a}z^{(1)}\widetilde{x}\widetilde{z}x^{(1)}x^{(0)},\phi)\psi_{\ell-p}(uz^{(1)}\widetilde{z})d(\cdots) .$$

Repeat the same steps, only that now, consider the function on $F^p\backslash\mathbb{A}^p$

(2.28)

$$b\mapsto \int\limits_{Z_F^{(1)}\backslash Z_{\mathbb{A}}^{(1)}} \int\limits_{F^p\backslash\mathbb{A}^p} \int\limits_{U_{\ell-p}(F)\backslash U_{\ell-p}(\mathbb{A})} E_1(u\widetilde{a}z^{(1)}\widetilde{\widetilde{b}}h,\phi)\psi_{\ell-p}(uz^{(1)})dudadz^{(1)}$$

for $h$ fixed and

$$\widetilde{\widetilde{b}} = \begin{pmatrix} I_p & 0 & b &&&& \\ & I_{2n-\ell-2} & 0 &&&& \\ && 1 &&&& \\ &&& I_{2(\ell-p)+3} &&& \\ &&&& 1 & 0 & b' \\ &&&&& I_{2n-l-2} & 0 \\ &&&&&& I_p \end{pmatrix} .$$

Write the Fourier expansion of the function (2.28), as in (2.26), and use

$$\eta^{(2)} = \begin{pmatrix} I_p &&&&&& \\ 0 & I_{2n-\ell-3} &&&&& \\ \eta & 0 & 1 &&&& \\ &&& I_{2(\ell-p)+5} &&& \\ &&&& 1 && \\ &&&& 0 & I_{2n-\ell-3} & \\ &&&& \eta' & 0 & 1 \end{pmatrix}$$

instead of $\eta^{(1)}$ in the previous step. Split $\widetilde{\mathcal{L}} = \mathcal{L}^{(2)} \cdot \widetilde{\widetilde{\mathcal{L}}}$, $\widetilde{Z} = Z^{(2)}\widetilde{\widetilde{Z}}$. $\mathcal{L}^{(2)}$ consists of all matrices (2.20) in $\widetilde{\mathcal{L}}$, such that $L$ has all lines zero, except for line $2n-\ell-2$, while $\widetilde{\widetilde{\mathcal{L}}}$ corresponds to $L$ with lines $2n-\ell-2, 2n-\ell-1$ and $2n-\ell$, being zero. $Z^{(2)} = \left\{ \begin{pmatrix} I_{2n-\ell-2} & * & 0 \\ & 1 & 0 \\ & & 1 \end{pmatrix} \right\}$, $\widetilde{\widetilde{Z}} = \left\{ \begin{pmatrix} \zeta & 0 & 0 \\ & 1 & 0 \\ & & 1 \end{pmatrix} \in Z_{2n-\ell} \right\}$.

Note that for $z^{(2)} \in Z^{(2)}$, $\eta^{(2)}\widetilde{\widetilde{b}}\eta^{(2)^{-1}} z^{(2)} = \widetilde{\widetilde{b}} z'^{(2)}$, for $z'^{(2)} \in Z^{(2)}$, such that $\psi(\eta \cdot b)\psi_{\ell-p}(z^{(2)}) = \psi_{\ell-p}(z'^{(2)})$, and so on. We repeat this procedure until we get

$$\int_{\mathcal{L}_{\mathbb{A}}} \int_{U^*_{\ell-p}(F)\backslash U^*_{\ell-p}(\mathbb{A})} E_1(vx, \phi)\psi_{\ell-p}(v)dvdx \tag{2.29}$$

where $U^*_{\ell-p}$ is the subgroup of matrices in $N_{\ell-p}$
(2.30)

$$\begin{pmatrix} I_p & A & X & D & Y \\ & Z & B & C & D' \\ & & I_{2(\ell-p)+1} & B' & X' \\ & & & Z^* & A' \\ & & & & I_p \end{pmatrix}, \quad \text{with } A \text{ having zero in the first column .}$$

Consider the function

$$\begin{pmatrix} m & e \\ 0 & 1 \end{pmatrix} \mapsto \int_{U^*_{\ell-p}(F)\backslash U^*_{\ell-p}(\mathbb{A})} E_1(v \begin{pmatrix} m & e \\ 0 & 1 \end{pmatrix}^{\wedge} h, \phi)\psi_{\ell-p}(v)dv \tag{2.31}$$

where $m \in \mathrm{GL}_p(\mathbb{A})$, $e \in \mathbb{A}^p$ (columns...) and

$$\begin{pmatrix} m & e \\ 0 & 1 \end{pmatrix}^{\wedge} = \begin{pmatrix} m & e & & & \\ & 1 & & & \\ & & I_{4n-2p-1} & & \\ & & & 1 & e' \\ & & & & m^* \end{pmatrix}.$$

The function (2.31) is "automorphic" and "cuspidal," in the sense that all constant terms are zero, as follows from Proposition 9. Thus its Fourier expansion (as for cusp forms of $\mathrm{GL}_{p+1}(\mathbb{A})$) yields that (2.29) is equal to

$$\sum_{\gamma \in Z_p(F)\backslash \mathrm{GL}_p(F)} \int_{\mathcal{L}_{\mathbb{A}}} \int_{N_{\ell-p}(F)\backslash N_{\ell-p}(\mathbb{A})} E_1(v\widehat{\gamma}x, \phi)\psi_{\ell-p}(v)dvdx \tag{2.32}$$

where $\widehat{\gamma} = \begin{pmatrix} \gamma & & \\ & I_{4n-2p+1} & \\ & & \gamma^* \end{pmatrix}$.

Let us record the identity that we proved

$$(2.33)\quad \int_{R_p(F)\backslash R_p(\mathbb{A})} \int_{N_\ell(F)\backslash N_\ell(\mathbb{A})} E_1(vr,\phi)\psi_\ell(v)dvdr$$
$$= \sum_{\gamma\in Z_p(F)\backslash \mathrm{GL}_p(F)} \int_{\mathcal{L}_{\mathbb{A}}} \int_{N_{\ell-p}(F)\backslash N_{\ell-p}(\mathbb{A})} E_1(v\widehat{\gamma}x\beta,\phi)\psi_{\ell-p}(v)dvdx$$

or, in the notation (2.11),

$$(2.34)\quad \int_{R_p(F)\backslash R_p(\mathbb{A})} p_\ell(r,\phi)dr = \sum_{\gamma\in Z_p(F)\backslash \mathrm{GL}_p(F)} \int_{\mathcal{L}_{\mathbb{A}}} p_{\ell-p}(\widehat{\gamma}x,\beta\cdot\phi)dx\ .$$

By our assumption $p_k(h,\phi)\equiv 0$, for $1\le k<\ell$. Note that $p_0(1,\phi)$ is the Whittaker coefficient of $E_1(\cdot,\phi)$, which is zero (Proposition 9). Thus, (2.34) implies that the constant terms of $p_\ell$ along the unipotent radicals $R_p$ are zero, for $1\le p<\ell$. This completes the proof in case $G=\mathrm{SO}_{4n+1}$. □

*Case* $G=\mathrm{Sp}_{4n}$. Let $R_p$ denote the unipotent radical of the standard, maximal parabolic subgroup of $\mathrm{Sp}_{2\ell}$, whose Levi part is isomorphic to $\mathrm{GL}_p\times \mathrm{Sp}_{2(\ell-p)}$ $(1\le p<\ell)$. We identify the elements of $R_p$, with their image in $\mathrm{Sp}_{4n}$, i.e. with matrices of the form

$$\begin{pmatrix} I_{2n-\ell} &&&& \\ & I_p & X & Y & \\ && I_{2(\ell-p)} & X' & \\ &&& I_p & \\ &&&& I_{2n-\ell}\end{pmatrix}.$$

We want to prove that

$$\int_{R_p(F)\backslash R_p(\mathbb{A})} p_\ell(r,\phi,\xi)dr\equiv 0,\quad \text{for all}\quad 1\le p<\ell\ ,$$

which, by (2.13) means that

$$(2.35)\quad \int_{R_p(F)\backslash R_p(\mathbb{A})} \int_{N_\ell(F)\backslash N_\ell(\mathbb{A})} E_1(vr,\phi)\theta^\xi_{\psi,\ell}(j_\ell(v)r)\psi_\ell(v)dvdr\equiv 0\ .$$

As in the previous case, it is convenient to conjugate $vr$ in (2.35) by the Weyl element

$$\beta=\begin{pmatrix} & I_p &&& \\ I_{2n-\ell} &&&& \\ && I_{2(\ell-p)} && \\ &&&& I_{2n-\ell} \\ &&& I_p & \end{pmatrix}.$$

We get matrices of the form

$$\beta v r \beta^{-1} = \begin{pmatrix} I_p & 0 & X & B & Y \\ L & Z & A & C & B' \\ & & I_{2(\ell-p)} & A' & X' \\ & & & Z^* & 0 \\ & & & L' & I_p \end{pmatrix} \tag{2.36}$$

where $z \in Z_{2n-\ell}$. Note that

$$j_\ell(v) = (L_{2n-\ell}, A_{2n-\ell}, B'_{2n-\ell}; C_{2n-\ell,1}) \tag{2.37}$$

$$\psi_\ell(v) = \psi(Z) = \psi(\sum_{j=1}^{2n-\ell-1} Z_{j,j+1}) \ .$$

Denote by $U'$ the group of matrices of the form (2.36), with $L, X, Y$ being zero and $Z = I_{2n-\ell}$. Denote by $\mathcal{L}$ the subgroup of matrices of the form (2.36) which are lower triangular. Finally, we identify $Z_{2n-\ell}$ with the subgroup of matrices of the form (2.36) where $A$, $B$, $C$, $X$, $Y$, and $L$ are all zero. Put, for $v \in N_\ell(\mathbb{A})$,

$$\widetilde{j}_\ell(\beta v \beta^{-1}) = j_\ell(v)$$

$$\widetilde{\psi}_\ell(\beta v \beta^{-1}) = \psi_\ell(v) \ .$$

The integral (2.35) becomes

$$\int\limits_{R_p(F)\backslash R_p(\mathbb{A})} \int\limits_{Z_{2n-\ell}(F)\backslash Z_{2n-\ell}(\mathbb{A})} \int\limits_{\mathcal{L}_F\backslash\mathcal{L}_{\mathbb{A}}} \int\limits_{U'_F\backslash U'_{\mathbb{A}}} \tag{2.38}$$

$$\cdot \ E_1(u'xz\beta r, \phi)\theta^\xi_{\psi,\ell}(\widetilde{j}_\ell(u(x)r)\widetilde{\psi}_\ell(z)du'dxdzdr.$$

By definition,

$$\theta^\xi_{\psi,\ell}(\widetilde{j}_\ell(u'x)r) = \sum_{\substack{\eta\in F^p \\ t\in F^{\ell-p}}} \omega_{\psi,\ell}(\widetilde{j}_\ell(u'x)r)\xi(\eta,t)$$

$$= \sum_{\substack{\eta\in F^p \\ t\in F^{\ell-p}}} \omega_{\psi_\ell}((\eta,0,0;0)\widetilde{j}_\ell(u'x)r)\xi(0,t) \ .$$

By the formulae of $\omega_{\psi,\ell}$, it is easy to see that for fixed $\eta, x$ and $r$, the function

$$u' \mapsto \sum_{t\in F^{\ell-p}} \omega_{\psi,\ell}((\eta,0,0;0)\widetilde{j}_\ell(u'x)r)\xi(0,t)$$

is left-$U'_F$ invariant. Now it is easy to see that

$$\int_{U'_F\backslash U'_{\mathbb{A}}} E_1(u'xz\beta r,\phi)\theta^{\xi}_{\psi,\ell}(\widetilde{j}_\ell(u'x)r)du'$$
$$=\int_{U'_F\backslash U'_{\mathbb{A}}} \sum_{\eta\in F^p} E_1(\eta^{(0)}u'xz\beta r,\phi)\sum_{t\in F^{\ell-p}}\omega_{\psi,\ell}((\eta,0,0;0)\widetilde{j}_\ell(u'x)r)\xi(0,t)du'$$
$$=\sum_{\eta\in F^p}\int_{U'_F\backslash U'_{\mathbb{A}}} E_1(u'\eta^{(0)}xz\beta r,\phi)\sum_{t\in F^{\ell-p}}\omega_{\psi,\ell}(j_\ell(u'\eta^{(0)}x)r)\xi(0,t)du'\ ,$$

where, as in the previous case, $\eta^{(0)}$ lies in $\mathcal{L}_F$ and is such that, in the notation of (2.36), $L=\begin{pmatrix}0\\ \vdots\\ 0\\ \eta\end{pmatrix}$. After interchanging summation and integration, we changed variable $u'\mapsto {\eta^{(0)}}^{-1}u'\eta^{(0)}$. Substituting this in (2.38), we obtained

$$(2.39)\quad \int_{R_p(F)\backslash R_P(\mathbb{A})}\ \int_{Z_{2n-\ell}(F)\backslash Z_{2n-\ell}(\mathbb{A})}\ \int_{\mathcal{L}'_F\backslash\mathcal{L}_{\mathbb{A}}}\ \int_{U'_F\backslash U'_{\mathbb{A}}}$$
$$\cdot\ E_1(u'xz\beta r,\phi)\sum_{t\in F^{\ell-p}}\omega_{\psi,\ell}(\widetilde{j}_\ell(u'x)r)\xi(0,t)\cdot\widetilde{\psi}_\ell(z)du'dxdzdr\ .$$

Here $\mathcal{L}'$ is the subgroup of matrices in $\mathcal{L}$ such that $L_{2n-\ell}=0$ (in the notation of (2.36). Let $\mathcal{L}^{(0)}$ be the subgroup of matrices in $\mathcal{L}$, such that all rows of $L$ are zero except for the last row. Note that for $r\in R_p(\mathbb{A})$ and $t\in F^{\ell-p}$,

$$\omega_{\psi,\ell}(r)\xi(0,t)=\xi(0,t)\ .$$

Thus, "conjugating $r$ back" in (2.39), we get

$$(2.40)$$
$$\int_{\mathcal{L}^{(0)}_{\mathbb{A}}}\ \int_{Z_{2n-\ell}(F)\backslash Z_{2n-\ell}(\mathbb{A})}\ \int_{\mathcal{L}'_F\backslash\mathcal{L}_{\mathbb{A}}}\ \int_{U_{\ell-p}(F)U_{\ell-p}(\mathbb{A})}$$
$$\cdot\ E_1(u'xza^{(0)}\beta,\phi)\sum_{t\in F^{\ell-p}}\omega_{\psi,\ell}(\widetilde{j}_\ell(ux^{(0)}))\xi(0,t)\cdot\widetilde{\psi}_\ell(z)du'dxdzda^{(0)}\ .$$

Take $\xi$ of the form $\xi_1\otimes\xi_2$, where $\xi_1\in S(\mathbb{A}^p)$ and $\xi_2\in S(\mathbb{A}^{\ell-p})$. Then, in (2.40),

$$\omega_{\psi,\ell}(\widetilde{j}_\ell(ua^{(0)}))\xi(0,t)=\xi_1(a)\omega_{\psi,\ell-p}(j_{\ell-p}(u))\xi_2(t)\ .$$

Thus, (2.40) becomes

$$(2.41)\quad \int_{\mathcal{L}_{\mathbb{A}}^{(0)}} \xi_1(a) \int_{Z_{2n-\ell}(F)\backslash Z_{2n-\ell}(\mathbb{A})} \int_{\mathcal{L}'_F\backslash\mathcal{L}'_{\mathbb{A}}} \int_{U_{\ell-p}(F)\backslash U_{\ell-p}(\mathbb{A})}$$

$$\cdot\ E_1(uxza^{(0)}\beta,\phi)\theta^{\xi_2}_{\psi,\ell-p}(j_{\ell-p}(u))\psi_{\ell-p}(uz)dudxdzda^{(0)}\ .$$

Now we are exactly at the same situation as in (2.24). We use a sequence of Fourier expansions on $E_1$, as in the previous case. Note that when we complete $uz$ to a matrix in $N_{\ell-p}$, by inserting elements in the zero, superdiagonal coordinates, then $j_{\ell-p}(u)$ remains unchanged. Therefore, we finally get

$$(2.42)\quad \sum_{\gamma\in Z_p(F)\backslash \mathrm{GL}_p(F)} \int_{\mathcal{L}_{\mathbb{A}}} \int_{N_{\ell-p}(F)\backslash N_{\ell-p}(\mathbb{A})}$$

$$\cdot\ E_1(v\widehat{\gamma}x\beta,\phi)\theta^{\xi_2}_{\psi,\ell-p}(j_{\ell-p}(v))\psi_{\ell-p}(v)\xi_1(j(x))dvdx\ ,$$

where for $x\in\mathcal{L}_{\mathbb{A}}$, we write, as before, $x=a^{(0)}\cdot u'$, with $u'\in\mathcal{L}'_{\mathbb{A}}$, and $j(x)=a$. Also, for $\gamma\in \mathrm{GL}_p(F)$, $\widehat{\gamma}=\begin{pmatrix}\gamma&&\\&I_{4n-2p}&\\&&\gamma^*\end{pmatrix}$. We record the identity that we proved (for $\xi=\xi_1\otimes\xi_2$)

$$(2.43)\quad \int_{R_p(F)\backslash R_p(\mathbb{A})} \int_{N_\ell(F)\backslash N_\ell(\mathbb{A})} E_1(vr,\phi)\theta^{\xi}_{\psi,\ell}(j_\ell(v)r)\psi_\ell(v)dvdr$$

$$= \sum_{\gamma\in Z_p(F)\backslash \mathrm{GL}_p(F)} \int_{\mathcal{L}_{\mathbb{A}}} \int_{N_{\ell-p}\backslash N_{\ell-p}(\mathbb{A})}$$

$$\cdot\ E_1(v\widehat{\gamma}x\beta,\phi)\theta^{\xi_2}_{\psi,\ell-p}(j_{\ell-p}(v))\psi_{\ell-p}(v)\xi_1(j(x))dvdx$$

or, in the notation (2.13),

$$(2.44)\quad \int_{R_p(F)\backslash R_p(\mathbb{A})} p_\ell(r,\phi,\xi_1\otimes\xi_2)dr = \sum_{\gamma\in Z_p(F)\backslash \mathrm{GL}_p(F)} \int_{\mathcal{L}_{\mathbb{A}}} p_{\ell-p}(\widehat{\gamma}x,\beta\phi,\xi_2)\xi_1(j(x))dx\ .$$

By our assumption $p_k(h,\phi,\xi)\equiv 0$, for $1\le k<\ell$. Note that in $p_0(1)=p_0(1,\phi)$ "there are no $\xi$ and no theta series." This is simply the Whittaker coefficient of $E_1(\cdot,\phi)$, which is zero (Proposition 9). Thus (2.44) implies that the constant terms of $p_\ell$ along the unipotent radicals $R_p$ are all zero. This completes the proof in case $G=\mathrm{Sp}_{4n}$. □

3. *Corollaries and conjectures.* We keep the notation of the previous sections. The first step of each of the towers (2.1), (2.2), (2.3), and (2.15) is clearly trivial, as follows from Proposition 9(2). In the first two towers, $\sigma_{2n-1}(\tau)$ (resp. $\sigma_{2n}(\tau)$) is just the restriction of $\overline{E_1(\cdot,\phi)}$ to $H_{2n-1}(\mathbb{A})$ (resp. $H_{2n}(\mathbb{A})$) and is hence nontrivial. In the last two towers, $\sigma_{2n-1}(\tau)$ (resp. $\sigma_{2n}(\tau)$) is easily seen to be nonzero for some choice of $\psi$. Thus, from Theorem 8, we conclude

THEOREM 10. *There is a natural number $\ell_n$, and a choice of $\psi$ in the last two cases, such that*

$$\ell_n \le \begin{cases} 2n-1, & \text{case (2.1)} \\ 2n, & \text{case (2.2)} \\ 2n-1, & \text{case (2.3)} \\ 2n, & \text{case (2.15)} \end{cases} , \tag{2.45}$$

$\sigma_k(\tau) = 0$, *for* $k < \ell_n$, *and* $\sigma_{\ell_n}(\tau)$ *is a nontrivial cuspidal representation of* $H_{\ell_n}(\mathbb{A})$.

We have the following(*)

CONJECTURE 1.

$$\ell_n = n \; ; \tag{2.46}$$

*i.e.*

$$H_{\ell_n} = \begin{cases} \mathrm{SO}_{2n+1}, & \text{case (2.1)} \\ \mathrm{SO}_{2n}, & \text{case (2.2)} \\ \widetilde{\mathrm{Sp}}_{2n}, & \text{case (2.3)} \\ \mathrm{Sp}_{2n}, & \text{case (2.15)} \end{cases} . \tag{2.47}$$

PROPOSITION 11. *$\sigma_{\ell_n}(\tau)$ is generic in the sense that it admits nontrivial Whittaker coefficients.*

*Proof.* Since $\sigma_{\ell_n}(\tau)$ is a nontrivial cuspidal representation, we can write it as a direct sum of irreducible, automorphic, cuspidal representations. Pick a summand $\sigma$. Then

$$\int\limits_{H_{\ell_n}(F)\backslash H_{\ell_n}(\mathbb{A})} \varphi(h)\overline{p_{\ell_n}(h)}dh \not\equiv 0 \tag{2.48}$$

(*) *Added in proof*. See the footnote, right after Theorem 8. Moreover, we can now also prove that $\sigma_k(\tau) = 0$ for $k < n$ in the remaining cases. This will appear in a work under preparation. Thus, Conjecture 1 is now proved.

for cusp forms $\varphi$, in the space of $\sigma$. As mentioned in (2.16), the integral (2.48) is the residue at $s = 1$ of a Rankin-Selberg integral for $H_{\ell_n} \times \mathrm{GL}_m$, where $m = 2n$, in cases (2.1)–(2.3) and $m = 2n+1$, in case (2.15). This integral represents (after a normalization) the standard partial $L$-function $L^S(\sigma\otimes\tau, s)$. (In case $H_k = \widetilde{\mathrm{Sp}}_{2k}$, the $L$-function depends on $\psi$.) However the Rankin-Selberg integral above is nonzero, if and only if $\sigma$ is generic. This is seen from its Euler product expansion (see [G-PS],[G], [So] and [G-R-S3]). In particular, it follows from (2.48) that $\sigma$ is generic. This implies that $\sigma_{\ell_n}(\tau)$ is generic. □

PROPOSITION 12. *Let $1 \le k < \ell_n$. Let $\sigma$ be an irreducible, automorphic, cuspidal, generic representation of $H_k(\mathbb{A})$. Then $L^S(\sigma\otimes\tau, s)$ is holomorphic at $s = 1$.*

*Proof.* By assumption, we have

$$\int\limits_{H_k(F)\backslash H_k(\mathbb{A})} \varphi(h)\overline{p_k(h)}dh \equiv 0 \tag{2.49}$$

for all cusp forms $\varphi$ in the space of $\sigma$. As in the proof of Proposition 11, the integral (2.49) is the residue at $s = 1$ of a Rankin-Selberg integral

$$\int\limits_{H_k(F)\backslash H_k(\mathbb{A})} \varphi(h)P_{k,s}(h)dh \tag{2.50}$$

which represents the standard partial $L$-function $L^S(\sigma\otimes\tau, s)$. In cases (2.1) and (2.2) (see (2.11)),

$$P_{k,s}(h) = P_{k,s}(h,\phi) = \int\limits_{N_k(F)\backslash N_k(\mathbb{A})} E(vh, f^{\phi}_{\tau,s})\psi_k(v)dv \tag{2.51}$$

and in cases (2.3) and (2.15) (see (2.13), (2.14)),

$$P_{k,s}(h) = P_{k,s}(h,\phi,\xi) = \int\limits_{N_k(F)\backslash N_k(\mathbb{A})} E(vh, f^{\phi}_{\tau,s})\theta^{\xi}_{\psi,k}(j_k(v)h)\psi_k(v)dv\ . \tag{2.52}$$

For decomposable data, the integral (2.50) equals an Euler product of the form

$$\prod_{\nu\in S} \mathcal{A}(W_\nu, P^{\nu}_{k,s}) \cdot L^S(\sigma\otimes\tau, s)\ , \tag{2.53}$$

where $S$ is a finite set of places, containing those at infinity, outside which the above data is unramified. $\mathcal{A}(W_\nu, P^{\nu}_{k,s})$ are the "local integrals" at the places of $S$. $W_\nu$ is a Whittaker function in a prescribed Whittaker model of $\sigma_\nu$, and $P^{\nu}_{k,s}$ is a corresponding local object obtained from $P_{k,s}$, depending on $\phi_\nu$ (resp. $(\phi_\nu, \xi_\nu)$). Thus the product (2.53) is holomorphic at $s = 1$ for all data at $S$.

We can choose data at $S$ (i.e. $W_\nu, \phi_\nu, \xi_\nu$) such that $\mathcal{A}(W_\nu, P^\nu_{k,s})$ is holomorphic and nonzero at a neighbourhood of $s = 1$ (see [So], [So1], [G-R-S3]).

We conclude that $L^S(\sigma\otimes\tau, s)$ is holomorphic at $s = 1$.(Recall, again, that for $H_k = \widetilde{\mathrm{Sp}}_{2k}$ this $L$-function depends on $\psi$.) □

We will prove in the next chapter one part of Conjecture 1, which is

THEOREM 13. *In case* (2.3),

$$\ell_n \geq n \ ;$$

*i.e. for an irreducible, automorphic, cuspidal representation* $\tau$ *of* $\mathrm{GL}_{2n}(\mathbb{A})$, *such that*

$$L^S(\tau, st, 1/2) \neq 0$$

*and*

$$L^S(\tau, \Lambda^2, s) \quad \textit{has a pole at } s = 1 \ ,$$

*we have*

$$\sigma_k(\tau) = 0, \quad \text{for} \quad 1 \leq k < n \ .$$

Combining this and Proposition 12, we get one of the main results of this paper.

THEOREM 14. *Let* $\sigma\otimes\tau$ *be an irreducible, automorphic, cuspidal, generic representation of* $\widetilde{\mathrm{Sp}}_{2k}(\mathbb{A}) \times \mathrm{GL}_{2n}(\mathbb{A})$, *such that* $1 \leq k < n$,

$$L^S(\tau, st, 1/2) \neq 0$$

*and*

$$L^S(\tau, \Lambda^2, s) \quad \textit{has a pole at} \quad s = 1 \ .$$

*Then* $L^S_\psi(\sigma\otimes\tau, s)$ *is holomorphic at* $s = 1$.

Finally, assume that $\ell_n$ is given by Conjecture 1. Let $\sigma$ be an irreducible summand of the cuspidal representation $\sigma_{\ell_n}(\tau)$. As in the proof of Proposition 11, and as in (2.53), keeping the same notation, we can find data at $S$, such that

$$\prod_{\nu\in S} \mathcal{A}(W_\nu, P^\nu_{\ell_n,s}) L^S(\sigma\otimes\tau, s)$$

has a pole at $s = 1$. We expect that, for all data, $\mathcal{A}(W_\nu, P^\nu_{\ell_n,s})$ is holomorphic at $s = 1$, so that $L^S(\sigma\otimes\tau, s)$ has a pole at $s = 1$. We cannot prove this at this stage. We make the following stronger conjecture.(*)

(*) *Added in proof*. See the footnotes after Theorems 8 and 10.

Conjecture 2. *The number $\ell_n$ is given by Conjecture* 1, *and $\sigma_{\ell_n}(\tau)$ is irreducible. Moreover, $\tau$ is the functorial lift of $\sigma_{\ell_n}(\tau)$ from $H_{\ell_n}$ to* $\mathrm{GL}_m$ (*where $m = 2n$ in cases* (2.1)–(2.3) *and $m = 2n+1$ in case* (2.15).)

*Remark* 1. Let $\sigma$ be an irreducible, automorphic, cuspidal generic representation of $H_{\ell_n}(\mathbb{A})$. Let $\tau$ be its Langlands functorial lift to $\mathrm{GL}_m(\mathbb{A})$ ($m$ as above). Assume that $\tau$ is cuspidal. Since $L^S(\sigma\otimes\tau, s)(= L^S(\tau\otimes\tau, s))$ has a pole at $s = 1$, we may apply our theory and find, by definition, that $\sigma$ is a summand of $\sigma_{\ell_n}(\tau)$. Thus $\sigma_{\ell_n}(\tau)$ contains all generic $\sigma$ which lift to $\tau$. Conjecture 2 implies that there is a unique generic $\sigma$ which lifts to $\tau$ (i.e. $\sigma = \sigma_{\ell_n}(\tau)$ is the unique generic member of its $L$-packet.)

*Remark* 2. Case (2.3) is a remarkable generalization of Waldspurger's result [W] on the theta correspondence between irreducible, automorphic, cuspidal representations $\tau$ of $\mathrm{PGL}_2(\mathbb{A})$ such that $L(\tau, \frac{1}{2}) \neq 0$ and $\psi$-generic representations of $\widetilde{\mathrm{SL}}_2(\mathbb{A})$. Our generalization is to $\tau$ – on $\mathrm{GL}_{2n}(\mathbb{A})$ such that $L^S(\tau, \Lambda^2, s)$ has a pole at $s = 1$ (this is automatic for $n = 1$ and $\tau$ with trivial central character) and $L(\tau, st, \frac{1}{2}) \neq 0$. To $\tau$ we associate $\sigma_n(\tau)$ on $\widetilde{\mathrm{Sp}}_{2n}(\mathbb{A})$.

## 3. Proof of vanishing of $\sigma_k(\tau)$, for $k < n$, in Case $G = \mathrm{Sp}_{4n}$

In this chapter, we prove Theorem 13. The main part of the proof is a local statement about the disjointness of $\mathrm{Sp}_{2n} \times \mathrm{Sp}_{2n}$-invariant functionals on representations of $G$ and certain eigenfunctionals with respect to a subgroup of $N_k$. For this chapter, we let $G = \mathrm{Sp}_{4n}$, $H = \mathrm{Sp}_{2n} \times \mathrm{Sp}_{2n}$, embedded in $G$ through $i$ in (1.7). We fix $1 \le k < n$.

1. *A preliminary reduction.* Let $\tau$ be an irreducible, automorphic, cuspidal (self-dual) representation of $\mathrm{GL}_{2n}(\mathbb{A})$, such that $L^S(\tau, st, \frac{1}{2}) \neq 0$ and $L^S(\tau, \Lambda^2, s)$ has a pole at $s = 1$. Consider the following subgroup, $N^{(k)}$, of $N_k$,

$$(3.1)\qquad N^{(k)} = \left\{ v = \begin{pmatrix} z & u & * & * & * \\ & 1 & 0 & y & * \\ & & I_{2k} & 0 & * \\ & & & 1 & u' \\ & & & & z^* \end{pmatrix} \in G \,\middle|\, z \in Z_{2n-(k+1)} \right\},$$

and its character

$$(3.2)\quad \chi_k(v) = \psi_k(v)\psi(y) = \psi(z_{12} + z_{23} + \cdots + z_{2n-k-2,2n-k-1} + u_{2n-k-1})\psi(y)\ .$$

Lemma 15. *We have $\sigma_k(\tau) = 0$, if and only if*

$$(3.3)\qquad \int\limits_{N^{(k)}(F)\backslash N^{(k)}(\mathbb{A})} E_1(v, \phi)\chi_k(v)dv \equiv 0\ .$$

*Proof.* The elements of $\sigma_k(\tau)$ are given by $p_k(h,\phi,\xi)$ in (2.13). $\sigma_k(\tau)=0$, if and only if $p_k(1,\phi,\xi)\equiv 0$ (for all $\phi$ and $\xi$). We have

$$\overline{p_k(1,\phi,\xi)} = \int\limits_{N_k(F)\backslash N_k(\mathbb{A})} E_1(v,\phi)\sum_{x\in F^k}\omega_{\psi,k}(j_k(v))\xi(x)\psi_k(v)dv \ . \tag{3.4}$$

Let

$$N_{k,1}=\left\{v=\begin{pmatrix} z & * & * & * & * & * \\ & 1 & e & 0 & y & * \\ & & I_k & 0 & 0 & * \\ & & & I_k & e' & * \\ & & & & 1 & * \\ & & & & & z^* \end{pmatrix}\in N_k\right\}, \tag{3.5}$$

$$N_{k,2}=\left\{a^*=\begin{pmatrix} I_{2n-k-1} & & & & & \\ & 1 & 0 & a & 0 & \\ & & I_k & 0 & a' & \\ & & & I_k & 0 & \\ & & & & 1 & \\ & & & & & I_{2n-k-1} \end{pmatrix}\in G\right\} . \tag{3.6}$$

Then $N_k=N_{k,1}\rtimes N_{k,2}$, and we can rewrite (3.4), using the notation in (3.5) and (3.6), as

$$\begin{aligned}
&\int\limits_{F^k\backslash\mathbb{A}^k}\ \int\limits_{N_{k,1}(F)\backslash N_{k,1}(\mathbb{A})} E_1(va^*,\phi)\sum_{x\in F^k}\omega_{\psi,k}(j_k(a^*))\xi(x+e)\psi_k(v)\psi(y)dvda \\
&\qquad = \int\limits_{F^k\backslash\mathbb{A}^k}\ \int\limits_{N^{(k)}(F)\backslash N_{k,1}(\mathbb{A})} E_1(va^*,\phi)\omega_{\psi,k}(j_k(a^*))\xi(e)\psi_k(v)\psi(y)dvda \\
&\qquad = \int\limits_{F^k\backslash\mathbb{A}^k}\ \int\limits_{N^{(k)}(F)\backslash N_{k,1}(\mathbb{A})} E_1(va^*,\phi)\psi(2ew_k{}^t a)\psi_k(v)\psi(y)\xi(e)dvda \ .
\end{aligned}$$

Here $w_k=\begin{pmatrix} & & 1 \\ & \cdot^{\cdot^{\cdot}} & \\ 1 & & \end{pmatrix}$. Consider $v$ in the last integral, and take it in the form (3.5). Write

$$v=v'\cdot\widetilde{e}$$

where $v' \in N^{(k)}$ and $\widetilde{e} = \begin{pmatrix} I_{2n-k-1} & & & & & \\ & 1 & e & & & \\ & & I_k & & & \\ & & & I_k & e' & \\ & & & & 1 & \\ & & & & & I_{2n-k-1} \end{pmatrix}$. We have

$$\widetilde{e}a^* = \begin{pmatrix} I_{2n-k-1} & & & & & \\ & 1 & 0 & 0 & 2ew_k{}^t a & \\ & & I_k & 0 & 0 & \\ & & & I_k & 0 & \\ & & & & 1 & \\ & & & & & I_{2n-k-1} \end{pmatrix} a^*\widetilde{e}\ . \tag{3.7}$$

Using this in the last integral, and changing variable $y \mapsto y - 2ew_k{}^t a$, we obtain

$$\overline{p_k(1,\phi,\xi)} = \int\limits_{e\in\mathbb{A}^k}\ \int\limits_{a\in F^k\backslash\mathbb{A}^k}\ \int\limits_{v'\in N^{(k)}(F)\backslash N^{(k)}(\mathbb{A})} E_1(v'a^*\widetilde{e},\phi)\chi_k(v')\xi(e)dv'dade\ .$$

Thus, for given $\phi$, $p_k(1,\phi,\xi) \equiv 0$, for all $\xi \in S(\mathbb{A}^k)$, if and only if

$$\int\limits_{F^k\backslash\mathbb{A}^k}\ \int\limits_{N^{(k)}(F)\backslash N^{(k)}(\mathbb{A})} E_1(v'a^*\widetilde{e},\phi)\chi_k(v')dv'da \equiv 0, \quad \text{for all} \quad e \in \mathbb{A}^k\ .$$

Thus, $p_k(1,\phi,\xi) \equiv 0$, for all $\phi$ and all $\xi$, if and only if

$$\int\limits_{F^k\backslash\mathbb{A}^k}\ \int\limits_{N^{(k)}(F)\backslash N^{(k)}(\mathbb{A})} E_1(va^*,\phi)\chi_k(v)dvda \equiv 0, \quad \text{for all } \phi\ . \tag{3.8}$$

Consider in (3.8) a right translation of $E_1$ by $\widetilde{e}$, where $e \in F^k$. Using (3.7), and changing variable $y \mapsto y + 2ew_k{}^t a$, and then $v \mapsto \widetilde{e}v\widetilde{e}^{-1}$, recalling that $E_1$ is left $G(F)$ invariant, we obtain from (3.8)

$$\int\limits_{F^k\backslash\mathbb{A}^k}\ \int\limits_{N^{(k)}(F)\backslash N^{(k)}(\mathbb{A})} E_1(va^*,\phi)\chi_k(v)\psi(2ew_k{}^t a)dvda \tag{3.9}$$

which is a general Fourier coefficient of the following function on $F^k\backslash\mathbb{A}^k$

$$a \mapsto \int\limits_{N^{(k)}(F)\backslash N^{(k)}(\mathbb{A})} E_1(va^*,\phi)\chi_k(v)dv\ . \tag{3.10}$$

Thus the integral (3.8) is identically zero (for all $\phi$), if and only if the integral (3.9) is identically zero, for all $\phi$ and all $e \in F^k$, which is equivalent to the function (3.10) being identically zero, for all $\phi$. This completes the proof of the lemma. □

Our aim now is to show (3.3), for $k < n$.

2. *$H$-invariant functionals and $(N^{(k)}, \chi_k)$-eigenfunctionals are disjoint.* We formulate the local result which is the heart of the proof of Theorem 13. In this section, $F$ is a local non-archimedean field.

THEOREM 16. *For $0 \leq k < n$, the Jacquet module $\mathcal{J}_{N^{(k)},\chi_k}(\mathrm{Ind}^{c\,G}_{\ H} 1)$ is zero.*

Note that $N^{(0)}$ is the standard maximal unipotent subgroup of $G$, and $\chi_0$ is its standard nondegenerate character. As a corollary, we get

THEOREM 17. *Let $\pi$ be an irreducible, admissible representation of $\mathrm{Sp}_{4n}(F)$. Assume that (the space of) $\pi$ admits nontrivial $H$-invariant functionals. Then, for $0 \leq k < n$, $\mathcal{J}_{N^{(k)},\chi_k}(\pi) = 0$, i.e. $\pi$ has no nontrivial $(N^{(k)}, \chi_k)$-eigenfunctionals.*

The reason is that $\pi$ admits nontrivial $H$-invariant functionals, if and only if $\widehat{\pi}$ admits nontrivial $H$-invariant functionals (soon to be explained) and this happens if and only if $\widehat{\pi}$ embeds in $\mathrm{Ind}_H^G 1 = C^\infty(H\backslash G)$, which is the same as the existence of a surjection

$$\mathrm{Ind}^{c\,G}_{\ H} 1 \longrightarrow \widehat{\widehat{\pi}} \cong \pi$$

and by the exactness of Jacquet functors and Theorem 16, we get $\mathcal{J}_{N^{(k)},\chi_k}(\pi) = 0$. Let us explain why $\pi$ admits nontrivial $H$-invariant functionals if and only if $\widehat{\pi}$ does. By [MVW, p. 91], $\widehat{\pi} \cong \pi^\delta$ where $\pi^\delta$ is the composition of $\pi$ with conjugation by $\delta^* = \begin{pmatrix} I_{2n} & \\ & \delta I_{2n} \end{pmatrix} \in G\mathrm{Sp}_{4n}(F)$, such that $\delta$ is not a square in $F$. It is clear that conjugation by $\delta^*$ preserves $H$ (see (1.7)). Thus, if $\pi$ acts on the space $V_\pi$ and $\ell$ is a $(\pi, H)$-invariant functional on $V_\pi$, then $\ell(\pi^\delta(h)\xi) = \ell(\pi(h^{\delta^*})\xi) = \ell(\xi)$ for $h \in H$ and $\xi \in V_\pi = V_{\pi^\delta}$, so that $\ell$ serves as an $H$-invariant functional for $\pi^\delta$ as well.

The proof of Theorem 16 will be given in Section 4 of this chapter.

3. *Theorem* 16 *implies Theorem* 13. We return to the global setup. In this section, $F$ is a number field and $\tau$ is as in Section 1. Consider the residue representation $E_1(\tau)$ in the space $\{E_1(g, \phi) = \operatorname*{Res}_{s=1} E(g, f^\phi_{\tau,s})\}$. It lies in the residual spectrum of $L^2(G(F)\backslash G(\mathbb{A}))$. Write $\overline{E_1(\tau)} = \oplus \pi_i$, a direct sum of irreducible representations (realized in concrete subspaces $V_{\pi_i}$). Let $k < n$. Denote the Fourier coefficient (3.3) by $W_k = W_k(\overline{E_1(\cdot, \phi)})$. We want to show that $W_k$ is zero on each summand $\pi_i$. By Theorem 2 and Corollary 3, the period along $H$ is nontrivial on $E_1(\tau)$. Therefore, there is a summand $\pi_{i_0}$, such that the period along $H$ is nontrivial on (the space of) $\pi_{i_0}$. Write $\pi_{i_0} \simeq \otimes \pi_{i_0,\nu}$, a tensor product of local representations. Pick a decomposable vector $\otimes f^0_\nu$, so that the period is nontrivial on its corresponding automorphic form in $V_{\pi_{i_0}}$. Let $S$ be a finite set of places such that for $\nu \notin S$, $f_\nu^0$ is the (prechosen) unramified vector of $\pi_{i_0,\nu}$. Fix any $\nu_0 \notin S$, and consider the restriction of the period along $H$ on the subspace of $V_{\pi_{i_0}}$, which corresponds to decomposable vectors $\otimes f_\nu$, where, for $\nu \neq \nu_0$, $f_\nu = f_\nu^0$, and $f_{\nu_0}$ varies in the space of $\pi_{i_0,\nu_0}$. This defines a nontrivial $H(F_{\nu_0})$-invariant functional on $\pi_{i_0,\nu_0}$, and hence, by Theorem 17, $\mathcal{J}_{N_{\nu_0}^{(k)},\chi_{k,\nu_0}}(\pi_{i_0,\nu_0}) = 0$. This implies that $W_k$ is trivial on $\pi_{i_0}$; otherwise, using the same reasoning as above, $W_k$ defines a nontrivial $(N_{\nu_0}^{(k)}, \chi_{k,\nu_0})$-eigenfunctional on $\pi_{i_0,\nu_0}$, that is a nontrivial element of $\mathcal{J}_{N_{\nu_0}^{(k)},\chi_{k,\nu_0}}(\pi_{i_0,\nu_0})$, which is a contradiction. To complete the proof of Theorem 13 (using Theorem 16), it remains to prove that $W_k$ is trivial on all summands $\pi_i$. Indeed, let $\pi_i$, $i \neq i_0$, be any other summand. It is clear that $\pi_i$ and $\pi_{i_0}$ are locally isomorphic at almost all places. (For almost all $\nu$, $\pi_{i,\nu}$ and $\pi_{i_0,\nu}$ are isomorphic to the unique unramified constituent of $\mathrm{Ind}_{P_\nu}^{G_\nu} \tau_\nu \otimes |\det\cdot|^{1/2}$.) What we just proved implies, for almost all $\nu$,

$$\mathcal{J}_{N_\nu^{(k)},\chi_{k,\nu}}(\pi_{i,\nu}) \cong \mathcal{J}_{N_\nu^{(k)},\chi_{k,\nu}}(\pi_{i_0,\nu}) = 0$$

and, in particular, $W_k$ is trivial on $\pi_i$.

It now remains to prove Theorem 16.

4. *Proof of Theorem* 16. We return to the notation of Section 2. $F$ is a local non-archimedean field. $G = \mathrm{Sp}_{4n}(F)$, $N^{(k)} = N^{(k)}(F)$, etc. ($0 \le k < n$ is fixed). Let $T$ denote the diagonal subgroup of $G$ and $W$, the Weyl group of $G$.

We prove Theorem 16 by standard Bruhat theory. For this, we have to describe the double cosets $H\backslash G/N^{(k)}$, and show that for all $g \in G$,

$$\chi_k\Big|_{g^{-1}Hg \cap N^{(k)}} \not\equiv 1 \ . \tag{3.11}$$

We start by realizing $H\backslash G$ inside $G$. Let

$$\varepsilon = \begin{pmatrix} I_n & & \\ & -I_{2n} & \\ & & I_n \end{pmatrix}$$

and define

$$\theta(g) = \varepsilon g \varepsilon^{-1} = \varepsilon g \varepsilon \ .$$

Note that $\theta$ is of order two, and that $\theta$ preserves root subgroups. The centralizer of $\varepsilon$ in $G$ is $H$, i.e. $\theta(g) = g$, if and only if $g \in H$. Let

$$Y = \{g^{-1}\theta(g) \mid g \in G\} \ .$$

Then

$$Y \cong H\backslash G \ ,$$

by

$$f : Hg \mapsto g^{-1}\theta(g) \ .$$

The natural right action of $G$ on $H\backslash G$, is translated, via $f$ to a twisted conjugation of $G$ on $Y$

$$y \cdot g = g^{-1}y\theta(g) \ .$$

As a first step, we show, after [Sp] (see also [J-R2, Lemma 2]).

PROPOSITION 18. *Each $y \in Y$ can be written in the form*

$$y = u^{-1}wa\theta(u) \ , \tag{3.12}$$

*where $u \in N^{(0)}$, $w \in W$, $a \in T$, such that*

$$\theta(wa) = (wa)^{-1} \ . \tag{3.13}$$

($W$ *is chosen to be generated by the simple reflection matrices.*)

*Proof.* Since $\theta(N^{(0)}) = N^{(0)}$, we can write the Bruhat decomposition of $y$ in the form

$$y = u_1^{-1}wa\theta(u_2) \ , \quad u_i \in N^{(0)}, w \in W, a \in T \ . \tag{3.14}$$

By definition of $Y$, its elements satisfy

$$\theta(y) = y^{-1} \ , \quad y \in Y \ . \tag{3.15}$$

Applying (3.15) to (3.14), we get

$$\theta(u_1^{-1})\theta(wa)u_2 = \theta(u_2^{-1})(wa)^{-1}u_1 \ . \tag{3.16}$$

By the uniqueness property of Bruhat decomposition, and the fact that

$$\theta\Big(N_G(T)/T\Big) = N_G(T)/T \ ,$$

we get, from (3.16)

$$\theta(wa) = (wa)^{-1} \ ,$$

which is (3.13). By definition of $\theta$, we have

$$\theta(w) = t \cdot w \ , \tag{3.17}$$

where $t$ is diagonal, with $\pm 1$ along the diagonal. From (3.13), we get

$$\theta(w)a = w^{-1}(wa^{-1}w^{-1}) \ ;$$

hence

$$\theta(w) = t' \cdot w^{-1} \ , \tag{3.18}$$

where $t'$ is diagonal with $\pm 1$ along the diagonal. We conclude that

$$w^2 = \delta \ , \tag{3.19}$$

where $\delta$ is diagonal, such that $\delta^2 = I$. (i.e. $w$ is of order two in $N_G(T)/T$). Note that (3.19) implies, that when we write $w$ as a permutation matrix in $\mathrm{GL}_{4n}(F)$, up to signs, then $w_{ij} \neq 0$ if and only if $w_{ji} \neq 0$ ($w_{ij} = \pm w_{ji} = \pm 1$). Note that since $w$ is symplectic, we have $w_{k,i} \neq 0$, if and only if $w_{4n+1-k,4n+1-i} \neq 0$.

Consider, for a Weyl element $\gamma$,

$$N_\gamma^+ = \{u \in N^{(0)} \mid \gamma u \gamma^{-1} \in N^{(0)}\}$$
$$N_\gamma^- = \{u \in N^{(0)} \mid \gamma u \gamma^{-1} \in N_-^{(0)}\} \ .$$

($N_-^{(0)}$ is the opposite to $N^{(0)}$.) From (3.19), it is clear that

$$N_w^+ = N_{w^{-1}}^+ \quad \text{and} \quad N_w^- = N_{w^{-1}}^- \ . \tag{3.20}$$

It is also clear that

$$\theta(N_w^\pm) = N_w^\pm \ . \tag{3.21}$$

This follows from (3.17), (3.18), and the definition of $\theta$. We may assume that $u_2$, in (3.14), lies in $N_{w^{-1}}^-$. Write, in (3.14), $u_1 = v_1^+ v_1^-$, where $v_1^+ \in N_{w^{-1}}^+$ and $v_1^- \in N_{w^{-1}}^-$. Writing (3.15) again,

$$\theta(v_1^-)^{-1}\theta(v_1^+)^{-1}\theta(wa)u_2 = \theta(u_2^{-1})(wa)^{-1}v_1^+v_1^- \ . \tag{3.22}$$

Since $(wa)^{-1}v_1^+(wa) \in a^{-1}w^{-1}N_{w^{-1}}^+ wa \in N^{(0)}$, we conclude, from (3.22) and the uniqueness in the Bruhat decomposition that $u_2 = v_1^-$. Thus

$$y = (v_1^-)^{-1}((v_1^+)^{-1} \cdot wa)\theta(v_1^-) \ .$$

This means that the $N^{(0)}$-orbit of $y$ is the same as that of $(v_1^+)^{-1} \cdot wa$, and so we may assume that $y = (v_1^+)^{-1} \cdot wa$. Using (3.15) again, we get

$$\theta(v_1^+)^{-1}\theta(wa) = (wa)^{-1} \cdot v_1^+ \ .$$

Applying $\theta$, and using (3.16), we get

$$(v_1^+)^{-1} = \quad (wa)\theta(v_1^+)(wa)^{-1} \ . \tag{3.23}$$

Consider the map

$$\lambda(u) = wa\theta(u)(wa)^{-1} \ .$$

It is clear, from (3.13), (3.20) and (3.21), that $\lambda$ defines an automorphism of order two of $N_{w^{-1}}^+$. By (3.23), $\lambda(v_1^+) = (v_1^+)^{-1}$. Since $N_{w^{-1}}^+$ is nilpotent, there is $v \in N_{w^{-1}}^+$, such that

$$v_1^+ = \lambda(v^{-1}) \cdot v$$

and so

$$\begin{aligned} y &= (v_1^+)^{-1} \cdot wa = v^{-1}\lambda(v) \cdot wa \\ &= v^{-1} wa\theta(v)(wa)^{-1} \cdot wa = v^{-1} \cdot wa \cdot \theta(v) \ . \end{aligned}$$

This proves Proposition 18. □

Let

$$U^{(k)} = \left\{ \begin{pmatrix} I_{2n-k-1} & & \\ & u & \\ & & I_{2n-k-1} \end{pmatrix} \middle| \begin{array}{l} u \text{ lies in the standard} \\ \text{maximal unipotent subgroup} \\ \text{of } \mathrm{Sp}_{2k+2}(F) \end{array} \right\} .$$

Clearly

$$N^{(0)} = U^{(k)} N^{(k)} ,$$

$U^{(k)}$ normalizes $N^{(k)}$ and for $u \in U^{(k)}$, $v \in N^{(k)}$, we have

$$\chi_k(uvu^{-1}) = \chi_k(v) \ . \tag{3.24}$$

Note that $U^{(k)} \cap N^{(k)}$ is the center of $U^{(k)}$. From Proposition 18, we have, of course, that the $N^{(k)}$-orbits in $Y$ have representatives of the form

$$u^{-1} \cdot wa \cdot \theta(u) , \tag{3.25}$$

where $u \in U^{(k)}$ and $wa$ as in (3.13).

Now, in order to prove (3.11), we have to show that $\chi_k$ is nontrivial on the stabilizer, in $N^{(k)}$, of elements of the form (3.25). An element $v \in N^{(k)}$ is in the stabilizer of (3.25) if and only if

$$v^{-1} \cdot (u^{-1} wa\theta(u))\theta(v) = u^{-1} wa\theta(u) \ ;$$

i.e.

$$(uvu^{-1})^{-1} wa\theta(uvu^{-1}) = wa \ .$$

By (3.24), it is enough to assume that $u = I$ in (3.25). Thus, we have to solve

$$v^{-1} \cdot wa \cdot \theta(v) = wa \quad , \quad v \in N^{(k)} \ .$$

We find it more convenient to replace $v$ by $\theta(v)$. Of course, this is permissible. The last equation transforms to

$$w(ava^{-1})w^{-1} = \varepsilon v \varepsilon \ . \tag{3.26}$$

Thus, if $v$ is in the stabilizer of $wa$, $w$ must "preserve the root structure" of $v$.

Let us call $wa$, satisfying (3.13), *nonrelevant*, if equation (3.26) admits solutions $v$ in $N^{(k)}$, such that $\chi_k(v) \neq 1$. Otherwise, call $wa$ *relevant*. So, in order to prove (3.11), we must show that all $wa$ satisfying (3.13) are nonrelevant, if $k < n$.

Let $L_k$ be the set of roots "inside $N^{(k)}$". For a root $\alpha$, and $t \in F$ denote by $x_\alpha(t)$ the corresponding standard unipotent matrix. Let $\Delta_k \subset L_k$ be the subset of roots $\alpha$, such that $\chi_k(x_\alpha(t)) \not\equiv 1$. The following two lemmas are useful.

LEMMA 19. *Let $\alpha \in \Delta_k$. Assume that $\beta = w(\alpha) \in L_k \setminus \Delta_k$. Then wa is nonrelevant.*

*Proof.* It is similar to [J-R2, Prop. 2]. Let

$$v = x_\alpha(t)x_\beta(s)x_{\alpha+\beta}(e) \tag{3.27}$$

where $x_{\alpha+\beta}(e) \equiv 1$, in case $\alpha+\beta$ is not root. We will show that there are $t, s, e$, such that $v$ solves (3.26) and $\chi_k(v) \neq 1$. We have

$$wx_\alpha(t)w^{-1} = x_\beta(ct) \quad , \quad c = \pm 1 \ . \tag{3.28}$$

By (3.19),

$$wx_\beta(s)w^{-1} = x_\alpha(c's) \quad , \quad c' = \pm 1 \ , \tag{3.29}$$

so that,

$$x_\alpha(cc't) = w^2 x_\alpha(t) w^{-2} = x_\alpha(\alpha(w^2)t) \ .$$

We conclude,

$$\alpha(w^2) = cc' \ . \tag{3.30}$$

Similarly, if $\alpha+\beta$ is a root, then

$$wx_{\alpha+\beta}(e)w^{-1} = x_{\alpha+\beta}(c''e) \ , \quad c'' = \pm 1 \ . \tag{3.31}$$

Using (in our case, since $\alpha$ and $\beta$ have the same length)

$$[x_\alpha(t), x_\beta(s)] = x_{\alpha+\beta}(dts) \ ,$$

where $d = 1, 2$ (depending on $\beta$), we get

$$\begin{aligned} x_{\alpha+\beta}(c''dts) = [x_\beta(ct), x_\alpha(c's)] = [x_\alpha(c's), x_\beta(ct)]^{-1} = x_{\alpha+\beta}(dc'cts)^{-1} \\ = x_{\alpha+\beta}(-c'cdts) \ . \end{aligned}$$

Thus

$$c'' = -cc' = -\alpha(w^2) \ . \tag{3.32}$$

Now substitute (3.27) in (3.26) and use (3.28), (3.29) and (3.31). We then have to solve

(3.33)

$$x_\beta(c\alpha(a)t)x_\alpha(c'\beta(a)s)x_{\alpha+\beta}(c''\alpha\beta(a)e) = x_\alpha(\alpha(\varepsilon)t)x_\beta(\beta(\varepsilon)s)x_{\alpha+\beta}(\alpha\beta(\varepsilon)e) \ .$$

Write the right-hand side of (3.33) as

$$x_\beta(\beta(\varepsilon)s)x_\alpha(\alpha(\varepsilon)t)x_{\alpha+\beta}(\alpha\beta(\varepsilon)(e+dts)) \ .$$

Thus, we must solve

$$c\alpha(a)t = \beta(\varepsilon)s \ , \tag{3.34}$$

$$c'\beta(a)s = \alpha(\varepsilon)t \ , \tag{3.35}$$

and if $\alpha + \beta$ is a root, also

$$c''\alpha\beta(a)e = \alpha\beta(\varepsilon)(e + dts) \ . \tag{3.36}$$

To solve (3.34) and (3.35), we must have

$$\frac{c\alpha(a)}{\beta(\varepsilon)} = \frac{\alpha(\varepsilon)}{c'\beta(a)} \ ;$$

i.e.

$$cc'\alpha\beta(a) = \alpha\beta(\varepsilon) \ . \tag{3.37}$$

We have,

$$\begin{aligned}\beta(a\varepsilon) &= w\alpha(a\varepsilon) = \alpha(wa\varepsilon w^{-1}) = \alpha(\varepsilon)\alpha(\varepsilon wa\varepsilon w^{-1}) \\ &= \alpha(\varepsilon)\alpha(\theta(wa)w^{-1}) = \alpha(\varepsilon)\alpha((wa)^{-1}w^{-1}) \\ &= \alpha(\varepsilon)\alpha(a^{-1}w^{-2}) = \alpha(\varepsilon a^{-1})\alpha(w^{-2}) = \alpha(\varepsilon a^{-1})cc' \ .\end{aligned}$$

We used (3.13), (3.19) and (3.30); (3.37) now follows (Recall that $c, c' = \pm 1$.) Now solve (3.36)

$$(c''\alpha\beta(a) - \alpha\beta(\varepsilon))e = \alpha\beta(\varepsilon)dts \ .$$

Thus, we only have to make sure that $c''\alpha\beta(a) - \alpha\beta(\varepsilon) \neq 0$. Indeed, by (3.32) and (3.37),

$$c''\alpha\beta(a) - \alpha\beta(\varepsilon) = -cc'\alpha\beta(a) - \alpha\beta(\varepsilon) = -2\alpha\beta(\varepsilon) \neq 0 \ .$$

Thus, choose, $s = \frac{c\alpha(a)}{\beta(\varepsilon)}t$, $e = -\frac{cd\alpha(a)t^2}{2\beta(\varepsilon)}$. Then $v$ in (3.27) solves (3.26) and

$$\chi_k(v) = \psi(t) \not\equiv 1 \ .$$

This proves Lemma 19. □

LEMMA 20. *For $w$ and $a$ satisfying* (3.13),

$$w_{i,4n+1-i} = 0 \ .$$

*That is, $w$ has zero along the second main diagonal.*

*Proof.* Let $i$ be such that $w_{i,4n+1-i} \neq 0$. We may assume that $i \leq 2n$ and $w_{i,4n+1-i} = 1$. Then, clearly, $w_{4n+1-i,i} = -1$. Let us concentrate on the $\mathrm{SL}_2$-subgroup of $G$ embedded in the coordinates $(i,i)$, $(i, 4n+1-i)$ , $(4n+1-i, i)$ and $(4n+1-i, 4n+1-i, i)$ and equate these coordinates only in equation (3.13), $\varepsilon w a \varepsilon = a^{-1} w^{-1}$. Write $a = \mathrm{diag}(a_1, \dots, a_{2n}, a_{2n}^{-1}, \dots, a_1^{-1})$. Since $a_{4n+1-i} = a_i^{-1}$ we get

$$\begin{pmatrix} & 1 \\ -1 & \end{pmatrix} \begin{pmatrix} a_i & \\ & a_i^{-1} \end{pmatrix} = \begin{pmatrix} a_i^{-1} & \\ & a_i \end{pmatrix} \begin{pmatrix} & -1 \\ 1 & \end{pmatrix}$$

(note that $\varepsilon$ has no effect here), which is impossible. □

To prove (3.11), we first show

LEMMA 21. *If $wa$ is relevant, then*

$$w_{ij} = 0 \quad \textit{for all} \quad 1 \le i, j \le 2n - k \ ;$$

*hence* (*since $w$ is a Weyl matrix*) *$w$ has the form*

$$w = \begin{pmatrix} \overbrace{\begin{matrix} 0 & \cdots & 0 \\ \vdots & & \vdots \\ 0 & \cdots & 0 \end{matrix}}^{2n-k} & * \\ * & \underbrace{\begin{matrix} 0 & \cdots & 0 \\ \vdots & & \vdots \\ 0 & \cdots & 0 \end{matrix}}_{2n-k} \end{pmatrix} .$$

*Proof.* Assume that $w_{11} = 1$. Then

$$w = \begin{pmatrix} 1 & & \\ & w' & \\ & & 1 \end{pmatrix}$$

where $w'$ is a Weyl matrix in $\mathrm{Sp}_{4n-2}(F)$. Consider the simple root $\alpha_1$, such that

$$x_{a_1}(t) = \begin{pmatrix} 1 & t & & & \\ & 1 & & & \\ & & I_{4n-4} & & \\ & & & 1 & -t \\ & & & & 1 \end{pmatrix} .$$

Then $\alpha_1 \in \Delta_k$ and $w(\alpha_1) \in L_k \backslash \Delta_k$, unless $w$ has the form

$$w = \begin{pmatrix} 1 & & & & \\ & 1 & & & \\ & & w'' & & \\ & & & 1 & \\ & & & & 1 \end{pmatrix} .$$

Consider the root $\alpha_2$, such that

$$x_{\alpha_2}(t) = \begin{pmatrix} 1 &&&&&& \\ & 1 & t &&&& \\ && 1 &&&& \\ &&& I_{4n-6} &&& \\ &&&& 1 & -t & \\ &&&&& 1 & \\ &&&&&& 1 \end{pmatrix} .$$

Then $\alpha_2 \in \Delta_k$ and $w(\alpha_2) \in L_k \backslash \Delta_k$, unless $w$ has the form

$$w = \begin{pmatrix} I_3 && \\ & w''' & \\ && I_3 \end{pmatrix} ,$$

etc. By Lemma 19, $w$ is nonrelevant, when we get to the $i^{\text{th}}$ simple root of $\mathrm{GL}_{2n-k}$, $\alpha_i$, and $w(\alpha_i) \in L_k \backslash \Delta_k$ (clearly, $w(\alpha_i) \in L_k$), $1 \le i < 2n - n$. Otherwise,

$$w = \begin{pmatrix} I_{2n-k-1} &&&& \\ & 1 &&& \\ && \widetilde{w} && \\ &&& 1 & \\ &&&& I_{2n-k-1} \end{pmatrix}$$

with $\widetilde{w}$ – a Weyl matrix in $\mathrm{Sp}_{2k}(F)$. For such $w$, equation (3.13) implies, for $a = \mathrm{diag}\ (a_1, \ldots, a_{2n}, a_{2n}^{-1}, \ldots, a_1^{-1})$,

$$a_i^2 = 1 \quad , \quad i = 1, \ldots, 2n-k \ .$$

Thus (3.26) is solved with $v = x_\eta(y) = \begin{pmatrix} I_{2n-k-1} &&&& \\ & 1 & 0 & y & \\ && I_{2k} & 0 & \\ &&& 1 & \\ &&&& I_{2n-k-1} \end{pmatrix}$.

Note that $\chi_k(v) = \psi(y) \not\equiv 1$. Thus, $w_{11} = 0$. Now, note that

$$w_{2n-k,j} = 0 \qquad\qquad j = 1, \ldots, 2n-k$$

(and hence,

$$w_{i,2n-k} = 0 \qquad\qquad i = 1, \ldots, 2n-k) \ .$$

Indeed, if, say, $w_{2n-k,j} = 1$ for $j \le 2n-k$. Then, by (3.19) (and the following remark) $w(\eta) \in L_k \backslash \Delta_k$, or $w(\eta) = \eta$. In the first case,

$$x_{w(\eta)}(y) = \begin{pmatrix} I_{j-1} &&&& \\ & 1 & 0 & \pm y & \\ && I_{n-2j} & 0 & \\ &&& 1 & \\ &&&& I_{j-1} \end{pmatrix} .$$

Using Lemma 19, we see that $wa$ is nonrelevant. If $w(\eta)=\eta$, use the argument above, which shows that $a^2_{2n-k}=1$ and hence $x_\eta(y)$ solves (3.26), so that $wa$ is again nonrelevant. Now, assume that $w_{j1}=1$, for $1<j<2n-k$. Then, $w_{1j}\neq 0$. For simplicity, we assume that $w_{1j}=1$ ($w_{j1},w_{1j}=\pm 1$ are allowed, in general). It is easy to see that $w(\alpha_j)\in L_k\backslash\Delta_k$, unless $w(\alpha_j)=\alpha_1$, which takes place if and only if $w_{j+1,2}\neq 0$. Again, we assume for simplicity that $w_{j+1,2}=w_{2,j+1}=1$ (recall (3.19)). Now consider $w(\alpha_{j+1})$. Again, $w(\alpha_{j+1})\in L_k\backslash\Delta_k$, unless $w(\alpha_{j+1})=\alpha_2$. We continue in this manner, and get that $wa$ is nonrelevant, if we find $j'$ such that $j+j'<2n-k-1$ and $w(\alpha_{j+j'})\in L_k\backslash\Delta_k$ (and then use Lemma 19). Otherwise, the upper left $(2n-k)\times(2n-k)$ block of $w$ has the form (up to signs)

$$\begin{array}{r}
\\ \\ \\ \\ \\ \\ \\ j\longrightarrow \\ \\ \\ \\ 2n-k\longrightarrow
\end{array}
\begin{pmatrix}
 & & & \overset{j}{\downarrow} & & & \overset{2n-k}{\downarrow} \\
0 & \cdots & 0 & 1 & & \cdots & 0 \\
 & & & & 1 & & \\
\vdots & & & & & \ddots & \\
 & & & & & 1 & 0 \\
0 & & & & & & 0 \\
1 & & & & & & \\
 & 1 & & & & & \\
 & & \ddots & & & & \vdots \\
 & & & 1 & & & \\
0 & \cdots & 0 & 0 & \cdots & & 0
\end{pmatrix}.$$

Now it is clear that $w(\alpha_{2n-k-1})\in L_k\backslash\Delta_k$, and by Lemma 19, $wa$ is nonrelevant. We continue in the same manner, row by row, in the upper left $(2n-k)\times(2n-k)$ block, and get the lemma. □

LEMMA 22. *If $wa$ is relevant, then*

$$w_{ij}=0\ ,\quad 1\le i\le 2n-k\ ,\quad 2n+k+1\le j\le 4n$$

*and hence (by Lemma 21 and (3.19)) $w$ has the form*

$$(3.38)\qquad w=\begin{pmatrix}
\overbrace{\begin{matrix}0 & \cdots & 0\\ \vdots & & \vdots \\ 0 & \cdots & 0\end{matrix}}^{2n-k} & \overbrace{\begin{matrix}0 & \cdots & 0\\ \vdots & & \vdots \\ 0 & \cdots & 0\end{matrix}}^{2n-k} \\
\begin{matrix}0 & \cdots & 0\\ \vdots & & \vdots \\ 0 & \cdots & 0\end{matrix} & \begin{matrix}0 & \cdots & 0\\ \vdots & & \vdots \\ 0 & \cdots & 0\end{matrix}
\end{pmatrix}
\begin{matrix}\left.\vphantom{\begin{matrix}0\\ \vdots\\ 0\end{matrix}}\right\}2n-k \\ \left.\vphantom{\begin{matrix}0\\ \vdots\\ 0\end{matrix}}\right\}2n-k\end{matrix}.$$

*Proof.* By Lemma 19, we know that $w_{1,4n} = 0$. Assume that $w_{1,4n+1-j} \neq 0$, for $2 \leq j \leq 2n-k$. Note that by (3.19), $w_{4n+1-j,1} \neq 0$. It is clear that $w(\alpha_{j-1}) \in L_k$, and $w(\alpha_{j-1}) \in L_k \backslash \Delta_k$, unless $w_{2,4n+2-j} \neq 0$ (and hence $w_{4n+2-j,2} \neq 0$) in which case $w(\alpha_{j-1}) = \alpha_1$. In the first case, it follows from Lemma 19 that $wa$ is nonrelevant. In the second case, we have $w(\alpha_{j-2}) \in L_k$, and $w(\alpha_{j-2}) \in L_k \backslash \Delta_k$, unless $w_{3,4n+3-j} \neq 0$ (and hence $w_{4n+3-j,3} \neq 0$), in which case $w(\alpha_{j-2}) = \alpha_2$, and so on. If at a certain stage we find a simple root $\alpha_{j-j'}$, such that $w(\alpha_{j-j'}) \in L_k \backslash \Delta_k$, then we are done, by Lemma 19. Otherwise the upper-right $(2n-k) \times (2n-k)$ corner of $w$ has the form

$$w = \begin{pmatrix} 0 \cdots \pm 1 & \overbrace{\quad\quad\quad \cdots}^{j} & 0 \\ & \pm 1 & \\ \vdots & \ddots & \vdots \\ & & \pm 1 \\ & & \vdots \\ 0 & \cdots & 0 \end{pmatrix} \Big\} j$$

and $j$ is even (Lemma 20). Consider the coordinates of $w$ in positions $\left(\frac{j}{2}, 4n - \frac{j}{2}\right), \left(\frac{j}{2}, 4n - \frac{j}{2} + 1\right), \left(\frac{j}{2} + 1, 4n - \frac{j}{2}\right), \left(\frac{j}{2} + 1, 4n - \frac{j}{2} + 1\right)$ and the corresponding $\mathrm{Sp}_4$-Weyl submatrix $w_j$ of $w$

$$w_j = \begin{pmatrix} 0 & 0 & b & 0 \\ 0 & 0 & 0 & c \\ -c^{-1} & 0 & 0 & 0 \\ 0 & -b^{-1} & 0 & 0 \end{pmatrix} \quad , \quad b, c = \pm 1 \ .$$

Put $a_{\frac{j}{2}} = d$, $a_{\frac{j}{2}+1} = e$. Then, in $\theta(w) = \varepsilon w \varepsilon$, $w_j$ either remains unchanged, or it is multiplied by $-1$. The second case takes place only if $j = 2n$, and then $k = 0$.

Denote the submatrix of $\theta(w)$, which comes in place of $w_j$ by

$$\begin{pmatrix} 0 & 0 & b' & 0 \\ 0 & 0 & 0 & c' \\ -c' & 0 & 0 & 0 \\ 0 & -b' & 0 & 0 \end{pmatrix} ((b'c') = \pm(b,c) = \pm(\pm 1, \pm 1)) \ .$$

Now examine condition (3.13) in the above coordinates only. Then, we must have

$$d = -bc'e \ .$$

Clearly $w(\alpha_{j/2}) = \alpha_{j/2}$. Now equation (3.26) is solved with $v = x_{\alpha_{j/2}}(t)$. Indeed, we have

$$wax_{\alpha_{\frac{j}{2}}}(t)a^{-1}w^{-1} = wx_{\alpha_{\frac{j}{2}}}(-bc't)w^{-1} = x_{\alpha_{\frac{j}{2}}}(cc't)$$

while

$$\varepsilon x_{\alpha_{\frac{j}{2}}}(t)\varepsilon = x_{\alpha_{\frac{j}{2}}}\left(\alpha_{\frac{j}{2}}(\varepsilon)t\right) .$$

We have $\alpha_{\frac{j}{2}}(\varepsilon) = -1$, if and only if $j = 2n$, and similarly $c' = -c$, if and only if $j = 2n$. Otherwise $\alpha_{\frac{j}{2}}(\varepsilon) = 1$ and $c' = c$. Since $c^2 = 1$, $x_{\alpha_{\frac{j}{2}}}(t)$ solves equation (3.26) and since $\chi_k(x_{\alpha_{\frac{j}{2}}}(t)) = \psi(t)$, we find that $wa$ is nonrelevant. We continue in this manner, row by row, in the upper right $(2n-k)\times(2n-k)$ block, and get the lemma. □

So far, we did not use the fact that $k < n$. To conclude the proof of (3.11), and hence of Theorem 13, we note that if $k < n$, the matrix $w$ of type (3.38) cannot be invertible, given that it is a Weyl matrix in $\mathrm{Sp}_{4n}(F)$.

Tel Aviv University, Tel Aviv, Israel (D.G.) and (D.S.)
*E-mail addresses*: ginzburg@math.tau.ac.il and soudry@math.tau.ac.il

The Ohio State University, Columbus, OH (S.R.)
*E-mail address*: haar@math.ohio-state.edu